\RequirePackage{ifpdf}
\ifpdf 
\documentclass[pdftex]{sigma}
\else
\documentclass{sigma}
\fi

\numberwithin{equation}{section}

\newcommand{\bH}{\mathbb{H}^{cc}}        
\newcommand{\sH}{\mathbb{H}^{-c}}      
\newcommand{\HO}{\mathbb{H}} 
\newcommand{\dahc}{\ddot{\mathfrak H}^{\mathfrak c}} 
\newcommand{\daha}{\ddot{\mathfrak{H}}} 
\newcommand{\aH}{\mathfrak{H}}
\newcommand{\sdaha}{\ddot{\mathfrak H}^-} 

\newcommand{\Z}{ \mathbb Z }
\newcommand{\Zn}{ \mathcal Z}
\newcommand{\be}{\beta}
\newcommand{\al}{\alpha}
\newcommand{\h}{\mathfrak{h}}
\newcommand{\e}{\mathfrak{e}}
\newcommand{\wtd}{\widetilde}
\newcommand{\td}{\tilde}
\newcommand{\C}{ \mathbb C }
\newcommand{\Cl}{ {\mathcal C} }

\newtheorem{Theorem}{Theorem}[section]
\newtheorem{Lemma}[Theorem]{Lemma}
\newtheorem{Proposition}[Theorem]{Proposition}
\newtheorem{Corollary}[Theorem]{Corollary}
\newtheorem{Remark}[Theorem]{Remark}
\newtheorem{Definition}[Theorem]{Definition}
\newtheorem{Example}[Theorem]{Example}

\numberwithin{table}{section}

\begin{document}

\allowdisplaybreaks

\renewcommand{\thefootnote}{$\star$}

\renewcommand{\PaperNumber}{012}

\FirstPageHeading

\ShortArticleName{Hecke--Clif\/ford Algebras and Spin Hecke Algebras IV:~Odd Double Af\/f\/ine Type}

\ArticleName{Hecke--Clif\/ford Algebras and Spin Hecke Algebras IV:\\ Odd Double Af\/f\/ine Type\footnote{This paper is a contribution to the Special
Issue on Dunkl Operators and Related Topics. The full collection
is available at
\href{http://www.emis.de/journals/SIGMA/Dunkl_operators.html}{http://www.emis.de/journals/SIGMA/Dunkl\_{}operators.html}}}

\Author{Ta KHONGSAP and Weiqiang WANG}
\AuthorNameForHeading{T.~Khongsap and W.~Wang}

\Address{Department of Mathematics, University of Virginia, Charlottesville, VA 22904, USA}

\Email{\href{mailto:tk7p@virginia.edu}{tk7p@virginia.edu}, \href{mailto:ww9c@virginia.edu}{ww9c@virginia.edu}}
\URLaddress{\url{http://people.virginia.edu/~tk7p/},
\url{http://www.math.virginia.edu/~ww9c/}}

\ArticleDates{Received October 15, 2008, in f\/inal form January 22,
2009; Published online January 28, 2009}

\Abstract{We introduce an odd double af\/f\/ine Hecke algebra (DaHa) generated by
a classical Weyl group $W$ and two skew-polynomial subalgebras of
anticommuting generators. This algebra is shown to be Morita
equivalent to another new DaHa which are generated by $W$ and two
polynomial-Clif\/ford subalgebras. There is yet a third algebra
containing a spin Weyl group algebra which is Morita
(super)equivalent to the above two algebras. We establish the PBW
properties and construct Verma-type representations via Dunkl
operators for these algebras.}


\Classification{20C08}

\section{Introduction}

\paragraph{1.1.}
The Dunkl operator  \cite{Dun}, which is an ingenious mixture of
dif\/ferential and ref\/lection ope\-ra\-tors, has found numerous
applications to orthogonal polynomials, representation theory,
noncommutative geometry, and so on in the past twenty years. To a
large extent, the Dunkl operators helped to motivate the def\/inition
of double af\/f\/ine Hecke algebras of Cherednik, which have played
important roles in several areas of mathematics. In recent years,
the representation theory of a degenerate version of the double
af\/f\/ine Hecke algebra (known as the rational Cherednik algebra or
Cherednik--Dunkl algebra) has been studied extensively (\cite{EG,
DO}; see the review paper of Rouquier~\cite{Rou} for extensive
references).

In \cite{W1}, the second author initiated a program of constructing
the so-called spin Hecke algebras associated to Weyl groups with
nontrivial $2$-cocycles, by introducing the spin af\/f\/ine Hecke
algebra as well as the rational and trignometric double af\/f\/ine Hecke
algebras associated to the spin symmetric group of I.~Schur
\cite{Sch}. Subsequently, in a series of papers \cite{Kh,KW1,KW2,W2}, the authors have extended the constructions of~\cite{W1} in
several dif\/ferent directions.

The construction of \cite{W1, KW2} provided two (super)algebras
$\dahc_W$ and $\sdaha_W$ associated to any classical Weyl group
$W$ which are Morita super-equivalent in the sense of~\cite{W2}.
These algebras admit the following PBW type properties:
\begin{gather*}
 \dahc_W  \cong  \C [\h^*] \otimes \Cl_{\h^*} \otimes \mathbb{C }W
\otimes\C [\h],  \qquad
 \sdaha_W  \cong   \Cl \{\h^*\} \otimes \C W^- \otimes\C [\h].
\end{gather*}
Here we denote by $\h$ the ref\/lection representation of $W$, by
$\C [\h^*]$ the polynomial algebra on $\h^*$, by $\Cl_{\h^*}$ a
Clif\/ford algebra, by $\Cl \{\h^*\}$ a skew-polynomial algebra with
anti-commuting generators, and by $\C W^-$ the spin Weyl group
algebra associated to the element $-1$ in the Schur multiplier
$H^2 (W,\C^*)$.

In contrast to the  rational Cherednik algebra (cf.~\cite{EG, Rou})
which admits a nontrivial automorphism group, the construction of
the algebras $\dahc_W$ is asymmetric as $\dahc_W$ contains as
subalgebras one polynomial algebra and one polynomial-Clif\/ford
subalgebras $\C [\h^*] \otimes \Cl_{\h^*}$ (the polynomial-Clif\/ford
algebra also appeared in the af\/f\/ine Hecke--Clif\/ford algebra of type
$A$ introduced by Nazarov \cite{Naz}). Moreover, in type $A$ case,
$\dahc_W$ contains Nazarov's algebra as a subalgebra, see~\cite{W1}.

\paragraph{1.2.}
In the present paper, we introduce three new algebras $\bH_W$,
$\sH_W$ and $\HO_W$, which are shown to be Morita (super)equivalent
to each other and to have PBW properties as follows:
\begin{gather*}
\bH_W  \cong  \C [\h] \otimes \Cl_{\h} \otimes \mathbb{C }W \otimes \Cl_{\h^*} \otimes\C [\h^*],
    \\
\sH_W  \cong  \Cl \{\h\} \otimes\mathbb{C }W^- \otimes
    \Cl_{\h^*} \otimes\C [\h^*],
    \\
 \HO_W  \cong  \Cl \{\h\} \otimes \C W  \otimes \Cl \{\h^*\}.
\end{gather*}

A novel feature here is that the algebra $\bH_W$ contains two
isomorphic copies of the polynomial-Clif\/ford subalgebra and there
is an automorphism of $\bH_W$ which switches these two copies.
Similar remark applies to the algebra $\HO_W$. We further show
that the odd DaHa $\HO_W$ of type $A$ contains the degenerate
af\/f\/ine algebra of Drinfeld and Lusztig as a subalgebra (see~\cite{EG} for a~similar phenomenon).

It turns out that the number of parameters in the algebras $\bH_W$,
$\sH_W$ and $\HO_W$ is equal to one plus the number of conjugacy
classes of ref\/lections in $W$, which is the same as for the
corresponding rational Cherednik algebras and dif\/fers by one from
the algebras introduced in~\mbox{\cite{W1, KW2}}. However, in contrast to
the usual Cherednik algebras, we show that each of the algebras
$\bH_W$, $\sH_W$ and $\HO_W$ contain large centers and are indeed
module-f\/inite over their respective centers.

\paragraph{1.3.} In Section~\ref{sec:spinWeylClifford} we present a f\/inite
dimensional version of the Morita (super) equivalence of the DaHa
mentioned above, and introduce the necessary concepts such as spin
Weyl group algebras and Clif\/ford algebras associated to the
ref\/lection representation $\h$.

The Schur multipliers $H^2(W,\C^*)$ for f\/inite Weyl groups $W$
were computed by Ihara and Yokonuma \cite{IY} (cf.\ Karpilovsky
\cite[Theorem~7.2.2]{Kar}). For example, $H^2(W_{B_n},\C^*)
=\Z_2\times \Z_2\times \Z_2$ for $n\geq 4$. Given any f\/inite Weyl
group $W$ (not necessarily classical) and any $2$-cocycle $\alpha
\in H^2(W,\C^*)$, we establish a superalgebra isomorphism (in two
versions $+$, $-$)
\[
\dot{\Phi}_{\pm}^{\al}: \ \Cl_{\h^*} \rtimes_{\pm} \;\C W^{-\al}
    \stackrel{\simeq}{\longrightarrow} \Cl_{\h^*} \otimes \C W^{\al}.
    \]
For the purpose of the rest of this paper, only the case when $W$
is classical and $\al = \pm 1$ is needed. The special case when
$\al =-1$ was established in \cite{KW1}, and this special case was
in turn a generalization of a theorem of Sergeev and Yamaguchi for
symmetric group.

We construct and study the algebras $\bH_W, \sH_W$ and $\HO_W$ in
the next three sections, i.e., in
Sections~\ref{sec:DaHa:doubleClifford},
 \ref{sec:sDaHCa}, and  \ref{sec:oddDaHa},
respectively. Among other results, we establish the PBW properties
as mentioned earlier and construct Verma-like representations of
the three algebras via Dunkl ope\-ra\-tors. Note in particular that a
representation for $\HO_W$ (see Theorems~\ref{ADunkl:eta},
\ref{BDunkl:eta}, and~\ref{DDunkl:eta}) is realized on the
skew-polynomial algebra with anti-commuting Dunkl operators.
Anti-commuting Dunkl operators f\/irst appeared in \cite{W1}, also
cf.~\cite{KW2}. In a very recent work \cite{BB}, Bazlov and
Berenstein introduced a notion of braided Cherednik algebra where
anti-commuting Dunkl operators also make a natural appearance.
After the second author communicated to them our construction of
$\HO_W$ for type $A$, they have also produced a similar algebra in
their second version (cf.\
\cite[Corollary~3.7]{BB}).

Finally, in the Appendix~\ref{sec:Appendix}, we collect
the proofs of several lemmas stated in
Section~\ref{sec:DaHa:doubleClifford} and
Section~\ref{sec:oddDaHa}.

\section[Schur multipliers of Weyl groups and Clifford algebras]{Schur multipliers of Weyl groups and Clif\/ford algebras} \label{sec:spinWeylClifford}

\subsection{A distinguished double cover}
\label{subsec:spinWeyl}

As in \cite{KW1, KW2}, we shall be concerned about a distinguished double covering
$\wtd{W}$ of $W$:
\[
1 \longrightarrow \Z_2 \longrightarrow \wtd{W} \longrightarrow W
\longrightarrow 1.
\]
We denote by $\Z_2 =\{1,z\},$ and by $\td{t}_i$ a f\/ixed preimage of
the generators $s_i$ of $W$ for each $i$. The group $\wtd{W}$ is
generated by $z, \td{t}_1,\ldots, \td{t}_n$ with relations
\[
z^2 =1, \qquad
 (\td{t}_{i}\td{t}_{j})^{m_{ij}} =
 \left\{
\begin{array}{rl}
1, & \text{if } m_{ij}=1,3, \\
z, & \text{if } m_{ij}=2,4,6.
\end{array}
\right.
\]

The quotient algebra $\C W^- :=\C \wtd{W} /\langle z+1\rangle$ of
$\C \wtd{W}$ by the ideal generated by $z+1$ is called the {\em
spin Weyl group algebra} associated to $W$. Denote by $t_i \in \C
W^-$ the image of $\td{t}_i$. It follows that $\C W^-$ is
isomorphic to the algebra generated by $t_i$, $1\le i\le n$, subject
to the relations
\begin{equation*}
(t_{i}t_{j})^{m_{ij}} = (-1)^{m_{ij}+1} \equiv \left\{
\begin{array}{rl}
1, & \text{if } m_{ij}=1,3, \\
-1, & \text{if } m_{ij}=2,4,6.
\end{array}
\right.
\end{equation*}
The algebra $\C W^-$ has a natural
superalgebra (i.e.\ $\Z_2$-graded) structure by letting each $t_i$
be odd.

\begin{Example}\rm  \label{present}
Let $W$ be the Weyl group of type $A_{n}$, $B_{n}$, or $D_{n}$, which
will be considered extensively in later sections. Then the spin
Weyl group algebra $\C W^-$ is generated by $t_1,\ldots, t_n$ with
relations listed in Table~\ref{table1}.

\begin{table}[h!]
\centering
\caption{The def\/ining relations of $\C W^-$.}\label{table1} \vspace{1mm}

\begin{tabular}{|l|l|}\hline \tsep{3pt}\bsep{3pt} Type of $W$ & Def\/ining Relations for $\C W^-$\\
\hline \tsep{3pt} $A_{n}$ & $t_{i}^{2}=1$,
$t_{i}t_{i+1}t_{i}=t_{i+1}t_{i}t_{i+1}$,\\ &
$(t_{i}t_{j})^{2}=-1\text{ if } |i-j| >1$\bsep{3pt} \\
\hline  \tsep{3pt} & $t_{1},\ldots,t_{n-1}$ satisfy the relations for $\C W^-_{A_{n-1}}$, \\
$B_{n}$ & $t_{n}^{2}=1$, $(t_{i}t_{n})^{2}=-1$ if $i\neq n-1,n$, \\
&
$(t_{n-1}t_{n})^{4}=-1$\bsep{3pt}\\
\hline \tsep{3pt} & $t_{1},\ldots,t_{n-1}$ satisfy the relations for
$\C W^-_{A_{n-1}}$,\\
$D_{n}$& $t_{n}^{2}=1$, $(t_{i}t_{n})^{2}=-1$ if $i\neq n-2, n$, \\&
$t_{n-2}t_{n}t_{n-2}=t_{n}t_{n-2}t_{n}$ \bsep{3pt}\\ \hline
\end{tabular}
\end{table}

\end{Example}

\subsection[Clifford algebra]{Clif\/ford algebra} \label{subsec:cliff}

Denote by $\h$ the ref\/lection representation of the Weyl group
$W$ (i.e.~a Cartan subalgebra of the corresponding complex Lie
algebra $\mathfrak g$). In the case of type $A_{n-1}$, we will
always choose to work with the Cartan subalgebra $\h =\C^n$ of
$gl_n$ instead of $sl_n$ in this paper.

Note that $\h$ carries a $W$-invariant nondegenerate bilinear form
$(-,-)$, which gives rise to an identif\/ication $\h^*\cong \h$ and
also a bilinear form on $\h^*$ which will be again denoted by
$(-,-)$. We identify $\h^*$ with a suitable subspace of $\C^N$ in
a standard fashion (cf.\ e.g.\ \cite[Table in~2.3]{KW1}). Then
describe the simple roots $\{\alpha_i\}$ for $\mathfrak g$ using a
standard orthonormal basis $\{e_i\}$ of $\C^N$. It follows that
$(\alpha_i, \alpha_j) =-2\cos (\pi /m_{ij})$.

Denote by $\Cl_{\h^*}$ the Clif\/ford algebra associated to $(\h^*,
(-,-))$, which is regarded as a subalgebra of the Clif\/ford algebra
$\Cl_N$ associated to $(\C^N,(-,-))$. We shall denote by $c_i$ the
generator in $\Cl_N$ corresponding to $\sqrt{2}\e_i$ and denote by
$\be_i$ the generator of $\Cl_{\h^*}$ corresponding to the simple root~$\al_i$ normalized with $\be_i^2=1$. In particular,
$\mathcal{C}_{N}$ is generated by $c_{1},\ldots,c_{N}$ subject to
the relations
\begin{gather*} 
c_{i}^{2} =1,\qquad c_{i}c_{j} =-c_{j}c_{i} \qquad \text{if} \quad i\neq j.
\end{gather*}
For example, we have
\[
\be_{i}=\frac{1}{\sqrt{2}}(c_{i}-c_{i+1}), \qquad
1\leq i\leq n-1
\]
and an additional one
\[
\be_{n}=  \left \{
  \begin{array}{ll}
   c_{n}  &\quad \text{if  } W =W_{B_{n}},  \\
   \frac{1}{\sqrt{2}}(c_{n-1}+c_{n})  & \quad \text{if  } W =W_{D_{n}}.
 \end{array}
    \right.
\]
Note that $N=n$ in the above three cases. For a complete list of
$\beta_i$ for each Weyl group $W$, we refer to \cite[Section~2]{KW1}
for details.

The action of $W$ on $\h$ and $\h^*$ preserves the bilinear form
$(-,-)$ and thus $W$ acts as automorphisms of the algebra $\Cl_{\h^*}$.
This gives rise to a semi-direct product $\Cl_{\h^*} \rtimes \C W$.
Moreover, the algebra $\Cl_{\h^*} \rtimes \C W$ naturally inherits the
superalgebra structure by letting elements in $W$ be even and each
$\be_i$ be odd.

\subsection{A superalgebra isomorphism} \label{subsec:finite_iso}

We recall the following result of Morris (the type $A$ case goes
back to Schur).
\begin{Proposition}[\cite{Mo,Sch}]  \label{th:morris}
Let $W$ be a f\/inite Weyl group. Then, there exists a surjective
superalgebra homomorphism
$\Omega: \C W^-  {\longrightarrow} \Cl_{\h^*}$ which sends $t_i$
to $\be_i$ for each $i$.
\end{Proposition}

Given two superalgebras $\mathcal{A}$ and $\mathcal{B}$, we view
the tensor product of superalgebras $\mathcal{A}$ $\otimes$
$\mathcal{B}$ as a~superalgebra with multiplication def\/ined by
\[
(a\otimes b)(a^{\prime}\otimes b^{\prime})
=(-1)^{|b||a^{\prime}|}(aa^{\prime }\otimes bb^{\prime})\qquad (a,a^{\prime}\in\mathcal{A},\ b,b^{\prime}\in\mathcal{B}),
\]
where $|b|$ denotes the $\Z_2$-degree of $b$, etc.

Now, let  $\Cl_n \rtimes_- \C W^-$ denote the algebra generated by
the subalgebras $\Cl_n$ and $\C W^-$ with the following additional
multiplication:
\[
    t_i c_j = -c_j^{s_i}t_i \qquad \forall \, i,j.
\]
Note that $\Cl_n \rtimes_-  \C W^-$ has a natural superalgebra
structure by setting each $c_i$ and $t_j$ to be odd for all
admissible $i$, $j$. We also endow a superalgebra structure on
$\Cl_{\h^*} \otimes \C W$ by declaring all elements of $W$ to be
even.

\begin{Theorem}\label{th:isofinite:Cl_n}
We have an isomorphism of superalgebras:
 \[
 \dot{\Phi}: \Cl_{\h^*} \rtimes_-   \C W^-
    \stackrel{\simeq}{\longrightarrow} \Cl_{\h^*} \otimes \C W
 \]
which extends the identity map on $\Cl_{\h^*}$ and sends each
$t_i$ to $\be_i s_i$. The inverse map $\dot{\Psi}$ is an extension
of the identity map on $\Cl_{\h^*}$ and sends each $s_i$ to $\be_i
t_i$.
\end{Theorem}

We f\/irst prepare a few lemmas.

\begin{Lemma} \label{lem:braidmatch}
We have $(\dot{\Phi}(t_i)\dot{\Phi}(t_j))^{m_{ij}} = (-1)^{m_{ij}+1}.$
\end{Lemma}
\begin{proof}
    Proposition~\ref{th:morris} says that $(t_it_j)^{m_{ij}} \!=\!
(\be_i\be_j)^{m_{ij}} \!=\! (-1)^{m_{ij}+1}.$ Also recall that \mbox{$(s_i
s_j)^{m_{ij}}\! =\!1$}. Then we have
\begin{gather*}
    (\dot{\Phi}(t_i)\dot{\Phi}(t_j))^{m_{ij}} = (\be_i s_i \be_j s_j)^{m_{ij}}
    = (\be_i \be_j)^{m_{ij}}( s_i s_j)^{m_{ij}} = (-1)^{m_{ij}+1}.\tag*{\qed}
\end{gather*}\renewcommand{\qed}{}
\end{proof}

\begin{Lemma} \label{lem:smash}
We have $\be_j \dot{\Phi}(t_i) = -\dot{\Phi}(t_i) \, \be_j^{s_i}$ for all $i$, $j$.
\end{Lemma}
\begin{proof}
Note that $(\be_i,\be_i) =2\be_i^2=2$, and hence
\[
\be_j \be_i =-\be_i\be_j +(\be_j,\be_i) =-\be_i\be_j
+\frac{2(\be_j,\be_i)}{(\be_i,\be_i)} \be_i^2 =-\be_i \be_j^{s_i}.
\]
Thus, we have
\begin{gather*}
\be_j   \dot{\Phi}(t_i) = \be_j \be_i s_i
 = -\be_i \be_j^{s_i} s_i
 = -\be_i s_i \be_j^{s_i}
 = -\dot{\Phi}(t_i)   \be_j^{s_i}.\tag*{\qed}
\end{gather*}\renewcommand{\qed}{}
\end{proof}

\begin{proof}[Proof of Theorem~\ref{th:isofinite:Cl_n}]
The algebra $\Cl_{\h^*} \rtimes_-  \C W^-$ is generated by $\be_i$
and $t_i$ for all $i$. Lemmas~\ref{lem:braidmatch} and
\ref{lem:smash} imply that $\dot{\Phi}$ is a (super) algebra
homomorphism. Clearly $\dot{\Phi}$ is surjective, and thus an
isomorphism by a dimension counting argument.

Clearly, $\dot{\Psi}$ and $\dot{\Phi}$ are inverses of each other.
\end{proof}

Let us denote by $\Cl_{\h^* \oplus \h}$ the Clif\/ford algebra
associated to $\big(\big(\h^*, (-,-)\big) \oplus
\big(\h,(-,-)\big)\big)$, and regard it as a subalgebra of the
Clif\/ford algebra $\Cl_{2N}$ associated to $\big(\big(\C^N,
(-,-)\big) \oplus \big((\C^N)^*,(-,-)\big)\big)$. We shall denote by
$e_i$ and $\nu_i$ the counterparts to $c_i$ and $\be_i$ via the
identif\/ication $\Cl_{\h^*} \cong \Cl_{\h}$.

By \cite[Theorem 2.4]{KW1}, there exists an isomorphism of
superalgebras
\begin{gather} \label{iso:KW}
\Phi: \ \Cl_{\h} \rtimes \C W  \rightarrow \Cl_{\h} \otimes \C W^-
\end{gather}
which extends the identity map on $\Cl_{\h}$ and sends each $s_i$
to $-\sqrt{-1}\nu_i t_i$. The isomorphism $\Phi$ was due to
Sergeev and Yamaguchi when $W$ is the symmetric group.

\begin{Theorem} \label{th:isofinite:Cl_2n}
We have an isomorphism of superalgebras:
\[
\ddot{\Phi}: \   \Cl_{\h^* \oplus \h} \rtimes \C W
\stackrel{\simeq}{\longrightarrow} \Cl_{\h^* \oplus \h} \otimes \C W
\]
which extends the identity map on $\Cl_{\h^* \oplus \h}$ and sends
each $s_i$ to $\sqrt{-1}\be_i\nu_i s_i.$ The inverse map
$\ddot{\Psi}$ is the extension of the identity map on $\Cl_{\h^*
\oplus \h}$ which sends each $ s_i$ to $\sqrt{-1}\be_i\nu_i s_i.$
\end{Theorem}
\begin{proof}
The isomorphisms $\dot{\Phi}$ in Theorem \ref{th:isofinite:Cl_n}
and $\Phi$ in \eqref{iso:KW} can be readily extended to the
following isomorphisms of superalgebras which restrict to the
identity map on $\Cl_{\h^* \oplus \h}$:
\begin{gather*}
 \Phi:  \  \Cl_{\h^* \oplus \h} \rtimes \C W \stackrel{\simeq}{\longrightarrow}
\Cl_{\h} \otimes \left(\Cl_{\h^*} \rtimes_- \C W^-\right),\;\\
 \dot{\Phi}: \   \Cl_{\h} \otimes \left(\Cl_{\h^*} \rtimes_- \C W^-
\right) \stackrel{\simeq}{\longrightarrow} \Cl_{\h^* \oplus \h}
\otimes \C W.
\end{gather*}
Observe that $\ddot{\Phi} = \dot{\Phi}\circ \Phi$, and so
$\ddot{\Phi}$ is an isomorphism.
\end{proof}

\subsection{The case of general 2-cocycles}
The materials of this subsection generalize the
Section~\ref{subsec:finite_iso} above and \cite[Section~2]{KW1};
however, they will not be used in subsequent sections.

The Schur multipliers $H^2(W,\C^*)$ for f\/inite Weyl groups $W$
were computed by Ihara and Yokonuma \cite{IY} (also cf.\
Karpilovsky \cite[Theorem~7.2.2]{Kar}). In all cases, we have
$H^2(W,\C^*) \cong \prod\limits_{j=1}^k \Z_2$ for suitable $k=0,1,2,3$.

Consider the following central extension of $W$ by $H^2(W,\C^*)$:
\[
1\longrightarrow H^2(W,\C^*) \longrightarrow \wtd{W}
\longrightarrow W \longrightarrow 1.
\]
We denote by $z_i$ the generator of the $i$th copy of $\Z_2$ in
$H^2(W,\C^*) \cong \prod\limits_{j=1}^k \Z_2$ and by $t_i$ a f\/ixed
preimage of the generator $s_i$ of W for each $i$. The group
$\wtd{W}$ is generated by $z_1,\ldots,z_k$, $t_1,\ldots, t_n$
subject to that $z_i$ is central of order $2$ for all $i$, and the
additional relations shown in Table~\ref{table2} below (cf.\ \cite[Table~7.1]{Kar}). In particular, the values of $k$ can be read of\/f from
Table~\ref{table2}.

\begin{table}[h!]
\centering
\caption{Central extensions $\wtd{W}$ of Weyl groups.}\label{table2} \vspace{1mm}

\begin{tabular}{|l|l|} \hline \tsep{4pt}\bsep{3pt} Type of $W$ &
  Generators/Relations for $\wtd{W}$\hidewidth \\ \hline
\tsep{3pt}&$t_i^2 =1$, $1 \le i \le n,$ \\
$A_n \ (n \ge 3)$  &$
t_i t_{i+1} t_i = t_{i+1} t_i t_{i+1}$, $1 \le i \le n-1$ \\
&$ t_i t_j = z_1 t_j t_i \text{ if } m_{ij} =2$ \bsep{3pt}\\
\hline
\tsep{3pt} $B_2$
& $t_1^2 = t_2^2 =1$, $(t_1 t_2)^2 = z_{1} (t_2 t_1)^2$\bsep{3pt}\\
\hline
\tsep{3pt}
 $B_3$
& $t_1^2 = t_2^2 = t_3^2 =1$, $t_1 t_2 t_1 = t_2 t_1 t_2,$\\
&$t_1 t_3 =z_1 t_3 t_1$, $(t_2 t_3)^2 = z_2 (t_3 t_2)^2$\bsep{3pt}\\
\hline
\tsep{3pt}
&$t_i^2 =1$, $1 \le i \le n$, $t_i t_{i+1} t_i = t_{i+1} t_i t_{i+1}$, $1 \le i \le n-2$\\
$B_n \  (n \ge 4)$
&$t_i t_j = z_1 t_j t_i$, $1\le i<j \le n-1$, $m_{ij} =2$\\
&$t_i t_n =z_2 t_n t_i$, $1\le i \le n-2$\\
&$(t_{n-1} t_n)^2 = z_3 (t_n t_{n-1})^2$\bsep{3pt}\\
\hline
\tsep{3pt} $D_4$
& $t_i^2 =1$, $1 \le i \le 4$, $t_i t_j t_i = t_j t_i t_j \text{ if } m_{ij} =3$\\
&$t_1 t_3 =z_1 t_3 t_1$, $t_1 t_4 =z_2 t_4 t_1$, $t_3 t_4 =z_3 t_4 t_3$ \bsep{3pt}\\
\hline
\tsep{3pt} &$t_i^2 =1$, $1 \le i \le n$, $t_i t_j t_i = t_j t_i t_j \text{ if } m_{ij} =3$\\
$D_n \ (n \ge 5)$
& $t_i t_j =z_1 t_j t_i$,  $1 \le i<j \le n$, $m_{ij} =2$, $i \neq n-1$\\
&$t_{n-1} t_n =z_2 t_n t_{n-1} $\bsep{3pt}\\
\hline
\tsep{3pt} $E_{n =6,7,8}$
& $t_i^2 =1$, $1 \le i \le n$, $t_i t_j t_i = t_j t_i t_j  \text{ if } m_{ij} =3$\\
&$t_i t_j =z_1 t_j t_i,  \text{ if } m_{ij} =2$\bsep{3pt}\\
\hline
\tsep{3pt} &$t_i^2 =1$, $1 \le i \le 4$, $t_i t_{i+1} t_i = t_{i+1} t_i t_{i+1}$  $(i =1,3)$ \\
 $F_4$  & $t_i t_j =z_1 t_j t_i$,  $1 \le i<j \le 4$, $m_{ij} =2, $ \\
 &$(t_2 t_3)^2 = z_2 (t_3 t_2)^2$ \bsep{3pt}\\
 \hline
\tsep{3pt}
$G_2   $
& $t_1^2 =t_2^2 =1$, $(t_1 t_2)^3 = z_1(t_2t_1)^3$ \bsep{3pt}\\
\hline
\end{tabular}
\end{table}

For $\al =(\al_i)_{i=1,\ldots,k} \in H^2(W,\C^*)$, the quotient
$\C W^{\al} := \C \wtd{W}/\langle z_i - \al_i, \forall\, i \rangle$
can be identif\/ied as the algebra generated by $t_1,\ldots, t_n$
subject to the relations:
\[
(t_i t_j)^{m_{ij}} =
 \left\{
\begin{array}{rl}
1, & \text{if } m_{ij}=1,3, \\
\al_{ij}, & \text{if } m_{ij}=2,4,6,
\end{array}
\right.
\]
where $\al_{ij} \in \{\pm 1\}$ is specif\/ied by $\al \in H^2(W,\C^*)$
as in Table~\ref{table2}.

Let $\Cl_{\h^*} \rtimes_{-} \C W^{-\al}$ denote the algebra
generated by subalgebras $\Cl_{\h^*}$ and $\C W^{-\al}$ with the
following additional multiplication:
\[
t_i^- \be_j = -\be_j^{s_i} t_i^- \qquad \forall\, i,j,
\]
where we have denoted by $t_i^-$ the generators of the subalgebra
$\C W^{-\al}$ of $\Cl_{\h^*} \rtimes_{-} \C W^{-\al}$, in order to
distinguish from the generators $t_i$ of $\C W^{\al}$ below. We
impose superalgebra structures on the algebras $\Cl_{\h^*}
\rtimes_{-} \C W^{-\al}$ and on $\Cl_{\h^*} \otimes \C W^{\al}$ by
letting $t_i^-$ be odd, $t_i$ be even, and $\beta_i$  be odd for
all $i$.

\begin{Theorem} \label{th:isofinite:cocycle}
Fix a $2$-cocycle $\al \in H^2(W,\C^*)$. We have an isomorphism of
superalgebras:
 \[
 \dot{\Phi}_{-}^{\al}: \ \Cl_{\h^*} \rtimes_{-} \;\C W^{-\al}
    \stackrel{\simeq}{\longrightarrow} \Cl_{\h^*} \otimes \C W^{\al}
\]
which extends the identity map on $\Cl_{\h^*}$ and sends $t_i^-$
to $\be_i t_i$ for each $i$. The inverse map $\dot{\Psi}_{-}^{\al}$ is
the extension of the identity map on $\Cl_{\h^*}$ and sends $t_i$
to $\be_i t_i^-$ for each $i$.
\end{Theorem}
\begin{proof}
By Lemma~\ref{lem:smash}, we have $\be_j \be_i = -\be_i
\be_j^{s_i}$. Recall that $t_i^-$ is odd and $t_i$ is even. So we
have  $\be_j\dot{\Phi}_{-}^{\al}(t_i^-) =
-\dot{\Phi}_{-}^{\al}(t_i^-) \, \be_j^{s_i}$ for all admissible
$i$, $j$. Moreover,
\begin{gather*}
(\dot{\Phi}_{-}^{\al}(t_i^-)\dot{\Phi}_{-}^{\al}(t_j^-))^{m_{ij}}
 =(\be_i t_i \be_j t_j)^{m_{ij}}=(\be_i \be_j)^{m_{ij}}(t_i
t_j)^{m_{ij}} = (-1)^{m_{ij}+1}(t_i t_j)^{m_{ij}}\\
\phantom{(\dot{\Phi}_{-}^{\al}(t_i^-)\dot{\Phi}_{-}^{\al}(t_j^-))^{m_{ij}}}{} =
         \begin{cases}
            1  &\text{if } m_{ij}=1,3,\\
            -\al_{ij} &\text{if } m_{ij}=2,4,6.
         \end{cases}
    \end{gather*}
Clearly,  $\dot{\Phi}_{-}^{\al}$ preserves the $\Z_2$-grading.
Hence, it follows that $\dot{\Phi}_{-}^{\al}$ is a surjective
superalgebra homomorphism, and thus an isomorphism by dimension
counting. It is clear that $\dot{\Psi}_{-}^{\al}$ is the inverse
of $\dot{\Phi}_{-}^{\al}$.
\end{proof}

Denote by $\Cl_{\h^*} \rtimes_{+} \C W^{-\al}$ the algebra generated
by subalgebras $\Cl_{\h^*}$ and $\C W^{-\al}$ with the following
additional multiplication:
\[
t_i^+ \be_j = \be_j^{s_i} t_i^+  \qquad \forall\, i,j,
\]
where we have denoted by $t_i^+$ the generators of the subalgebra
$\C W^{-\al}$ of $\Cl_{\h^*} \rtimes_{+} \C W^{-\al}$, in order to
distinguish from the generators $t_i$ of $\C W^{\al}$. We impose
superalgebra structures on the algebras $\Cl_{\h^*} \rtimes_{+} \C
W^{-\al}$ and on $\Cl_{\h^*} \otimes \C W^{\al}$ by letting
$t_i^+$ be even, $t_i$ be odd, and $\beta_i$  be odd for all $i$.

\begin{Theorem}
Fix a $2$-cocycle $\al \in H^2(W,\C^*)$. We have an isomorphism of
superalgebras:
\[
\dot{\Phi}_{+}^{\al}: \ \Cl_{\h^*} \rtimes_{+} \;\C W^{-\al}
    \stackrel{\simeq}{\longrightarrow} \Cl_{\h^*} \otimes \C W^{\al}
    \]
which extends the identity map on $\Cl_{\h^*}$ and sends $t_i^+
\mapsto -\sqrt{-1}\be_i t_i$. The inverse map $\dot{\Psi}_{+}^{\al}$ is
the extension of the identity map on $\Cl_{\h^*}$ and sends each
$t_i$ to $\sqrt{-1}\be_i t_i^+$.
\end{Theorem}

\begin{proof}
By Lemma~\ref{lem:smash}, we have $\be_j \be_i = -\be_i
\be_j^{s_i}$. Recall that $t_i^+$ is even while $t_i$ is odd. Then
$\be_j\dot{\Phi}_{+}^{\al}(t_i^+) = \dot{\Phi}_{+}^{\al}(t_i^+) \,
\be_j^{s_i}$ for all admissible $i$, $j$. Moreover,
    \begin{gather*}
(\dot{\Phi}_{+}^{\al}(t_i^+)\dot{\Phi}_{+}^{\al}(t_j^+))^{m_{ij}}
 =(-\be_i t_i \be_j t_j)^{m_{ij}}= (\be_i \be_j t_i  t_j)^{m_{ij}}
          =(\be_i \be_j)^{m_{ij}}(t_i t_j)^{m_{ij}}\\
 \phantom{(\dot{\Phi}_{+}^{\al}(t_i^+)\dot{\Phi}_{+}^{\al}(t_j^+))^{m_{ij}}}{}          =
         \begin{cases}
            1  &\text{if } m_{ij}=1,3,\\
            -\al_{ij} &\text{if } m_{ij}=2,4,6.
         \end{cases}
    \end{gather*}
It follows that $\dot{\Phi}_{+}^{\al}$ is an isomorphism of
superalgebra with inverse $\dot{\Psi}_{+}^{\al}$.
\end{proof}

Denote by $\Cl_{\h^* \oplus \h} \rtimes_+ \C W^\al$ the algebra
generated by subalgebras $\Cl_{\h^*\oplus \h}$ and $\C W^{\al}$
with the following additional multiplication:
\[
t_i \be_j = \be_j^{s_i} t_i, \qquad t_i \nu_j = \nu_j^{s_i} t_i,
\qquad \forall\, i,j.
\]
We impose superalgebra structures on the algebras $\Cl_{\h^*
\oplus \h} \rtimes_+ \C W^{\al}$ and on $\Cl_{\h^* \oplus \h}
\otimes \C W^{\al}$ by letting each $t_i$ be even, and letting
each $\beta_i$, $\nu_i$ be odd.

\begin{Corollary}
For a $2$-cocycle $\al \in H^2(W,\C^*)$, we have an isomorphism of
superalgebras:
\[
\Cl_{\h^* \oplus \h} \rtimes_+ \C W^\al \cong  \Cl_{\h^* \oplus
\h} \otimes \C W^\al
\]
which extends the identity map on $\Cl_{\h^* \oplus \h}$ and sends
each $t_i$ to $\sqrt{-1} \beta_i\nu_i t_i$.
\end{Corollary}

\begin{Remark}\rm
When $\al = \pm1$ and $\C W^\al$ becomes the usual group
algebra $\C W$ or the spin group algebra $\C W^-$, we recover the
main results of Section~\ref{subsec:finite_iso}.
\end{Remark}

%
%
%
%
%

\section[The DaHa  with two polynomial-Clifford subalgebras]{The DaHa  with two polynomial-Clif\/ford subalgebras} \label{sec:DaHa:doubleClifford}

In the remainder of the paper, $W$ is always assumed to be one of
the classical Weyl groups of type $A_{n-1}$, $B_n$, or $D_n$, and we
shall often write $\Cl_{2n}$ for $\Cl_{\h^*\oplus \h}$.
\subsection[The definition of $\bH_W$]{The def\/inition of $\boldsymbol{\bH_W}$}

Let $W$ be one of the classical Weyl groups. The goal of this
section is to introduce a rational double af\/f\/ine Hecke algebra
(DaHa) $\bH_W$ which is generated by $\C W$ and two isomorphic
``polynomial-Clif\/ford'' subalgebras. Note that this construction is
dif\/ferent from the double af\/f\/ine Hecke--Clif\/ford algebra introduced
in \cite{W1, KW2} which is generated by $\C W$, a polynomial
subalgebra, and a ``polynomial-Clif\/ford'' subalgebra.

Identify $\C[\h^*] \cong \C[x_1,\ldots,x_n]$ and $\C[\h] \cong
\C[y_1,\ldots,y_n]$, where $x_i$, $y_i$ $(1\le i \le n)$  correspond
to the standard orthonormal basis $\{\e_i\}$ for $\h^*$ and its
dual basis $\{\e_i^*\}$ for $\h$, respectively. For~$x$,~$y$ in an
algebra $A$, we denote as usual that
\[
[x,y] =xy-yx \in A.
\]

\subsubsection[The algebra $\bH_W$ of type $A_{n-1}$]{The algebra $\boldsymbol{\bH_W}$ of type $\boldsymbol{A_{n-1}}$}

\begin{Definition}\rm
Let $t, u\in \C$ and $W=S_n$. The algebra $\bH_W$ of type
$A_{n-1}$ is generated by $x_i$, $y_i$  $(1\le i \le n)$, $\Cl_{2n}$
and $W$, subject to the following additional relations:
\begin{gather}
x_i x_j  =x_j x_i, \qquad y_i y_j =y_j y_i  \qquad (\forall\,  i,j),\nonumber \\
 \sigma c_i  =c_i^{\sigma} \sigma, \qquad \sigma e_i =e_i^{\sigma}
\sigma, \nonumber\\
\sigma x_i  = x_i^{\sigma} \sigma, \qquad \sigma y_i =
y_i^{\sigma} \sigma \qquad (\forall\,  \sigma \in W,\ \forall\, i),\label{commonA}\\
e_i x_j  = x_j e_i, \qquad c_i x_j = (-1)^{\delta_{ij}}x_j c_i \qquad (\forall\, i,j), \nonumber\\
c_i y_j  = y_j c_j, \qquad e_i y_j = (-1)^{\delta_{ij}} y_j e_i
\qquad (\forall\, i,j),\nonumber
\\
\lbrack y_i, x_j]  =  u (1 +c_ic_j)(1 +e_je_i)s_{ij} \qquad (i \neq j), \label{Ayixj}\tag{\rm 3.2a}\\
\lbrack y_i, x_i]  =  t c_ie_i
- u\sum_{k \ne i}  (1 +c_kc_i) (1 +e_ke_i)s_{ki}.\label{Ayixi}\tag{\rm 3.2b}
\end{gather}
\end{Definition}

Alternatively, we may view $t$, $u$ as formal variables and $\bH_W$
as a $\C[t,u]$-algebra. Similar remarks apply to other algebras
def\/ined in this paper.

\setcounter{equation}{2}

\subsubsection[The algebra $\bH_W$ of type $D_{n}$]{The algebra $\boldsymbol{\bH_W}$ of type $\boldsymbol{D_{n}}$}

Let $W =W_{D_n}$. Regarding elements in $W$ as even signed
permutations of $1,2,\ldots,n$ as usual, we identify the
generators $s_i \in W$, $1\leq i\leq n-1$, with transposition
$(i,i+1)$, and $s_n \in W$ with the transposition of $(n-1,n)$
coupled with the sign changes at $n-1$, $n$. For $1\le i\neq j\le n$,
we denote by $s_{ij} \equiv (i,j) \in W$ the transposition of $i$
and $j$, and $\overline{s}_{i j} \equiv \overline{(i,j)} \in W$
the transposition of $i$ and $j$ coupled with the sign changes at
$i$, $j$. By convention, we have
\[
\overline{s}_{n-1,n}\equiv \overline{(n-1,n)}=s_n,\qquad
\overline{s}_{ij}\equiv \overline{(i,j)}
=s_{jn}s_{i,n-1}s_{n}s_{i,n-1}s_{jn}.
\]

\begin{Definition}\rm
Let $t,u\in\C$ and $W =W_{D_n}$. The algebra $\bH_W$ of type
$D_{n}$ is generated by $x_i$, $y_i$  $(1\le i \le n)$, $\Cl_{2n}$
and $W$, subject to the relations (\ref{commonA}) with the current
$W$, and (\ref{Dyixj})--(\ref{Dyixi}) with $i\neq j$ below:
\begin{subequations}\label{Dyx}
\begin{gather}
\lbrack y_{i},x_{j}]  =  u (1 +c_ic_j)(1 +e_je_i)s_{ij}
-u(1-c_ic_j)(1-e_je_i)\overline{s}_{ij}, \label{Dyixj}  \\
\lbrack y_{i},x_{i}]   = tc_i e_i -u\sum_{k\neq
i}(1+c_{k}c_{i})(1+e_k e_i)s_{ki}
 +u\sum_{k\neq
i}(1-c_{k}c_{i})(1-e_ke_i)\overline{s}_{ki}. \label{Dyixi}
\end{gather}
\end{subequations}
\end{Definition}

\subsubsection[The algebra $\bH_W$ of type $B_{n}$]{The algebra $\boldsymbol{\bH_W}$ of type $\boldsymbol{B_{n}}$}

Let $W =W_{B_n}$. We identify $W$ as usual with the signed
permutations on $1, \ldots, n$. Regarding $W_{D_n}$ as a subgroup of
$W$, we have $s_{ij}, \overline{s}_{ij} \in W$ for $1 \le i\neq j
\le n$. Further denote $\tau_{i} \equiv \overline{(i)} \in W$ the
sign change at $i$ for $1\le i \le n$. By def\/inition, we have
\[
\tau_n \equiv \overline{(n)} =s_n, \qquad \tau_{i} \equiv
\overline{(i)} =s_{in}s_{n}s_{in}.
\]

\begin{Definition}\rm
Let $t,u,v\in\C$, and $W = W_{B_n}$. The algebra $\bH_W$ of type
$B_{n}$ is generated by $x_i$, $y_i$ $(1\le i \le n)$, $\Cl_{2n}$ and
$W$, subject to the relations (\ref{commonA}) with the current $W$,
and (\ref{Byixj})--(\ref{Byixi}) with $i\neq j$ below:
\begin{subequations}\label{Byx}
\begin{gather}
\lbrack y_{i},x_{j}]   =  u (1 +c_ic_j)(1 +e_je_i)s_{ij}
-u(1-c_ic_j)(1-e_je_i)\overline{s}_{ij},  \label{Byixj}  \\
\lbrack y_{i},x_{i}]    =  tc_i e_i -u\sum_{k\neq
i}(1+c_{k}c_{i})(1+e_k e_i)s_{ki}
 +u\sum_{k\neq
i}(1-c_{k}c_{i})(1-e_ke_i)\overline{s}_{ki} -v\tau_i. \label{Byixi}
\end{gather}
\end{subequations}
\end{Definition}

\subsection[The PBW basis for $\bH_W$]{The PBW basis for $\boldsymbol{\bH_W}$}

We shall denote $x^{\al}\!=\!x_1^{a_1}\cdots x_n^{a_n}$ for $\al\! =\!
(a_1,\ldots, a_n) \!\in\! \Z_{+}^n$, $c^\epsilon
\!=\!c_1^{\epsilon_1}\cdots c_n^{\epsilon_n}$ for \mbox{$\epsilon =
(\epsilon_1,\ldots,\epsilon_n) \!\in\! \Z_2^n$}. Similarly, we def\/ine
$y^\al$ and  $e^{\epsilon}$. Note that the algebra $\bH_W$
contains $\C[\h^*]$, $\Cl_{\h^*}$, $\C[\h]$, $\Cl_{\h},$ and $\C W$ as
subalgebras.

\begin{Theorem}  \label{PBW:ADBdahc}
Let $W$ be $W_{A_{n-1}}$, $W_{D_n}$ or $W_{B_n}$. The
multiplication of the subalgebras $\C[\h^*]$, $\C[\h]$,
$\Cl_{\h^*}$, $\Cl_{\h}$, and $\C W$ induces a vector space
isomorphism
\[
\C[\h^*] \otimes \Cl_{\h^*}  \otimes\C W \otimes \C[\h]  \otimes
\Cl_{\h} \stackrel{\simeq}{\longrightarrow} \bH_W.
\]
Equivalently, the elements $\{x^\al c^\epsilon  w e^{\epsilon'}
y^\gamma\,|\, \al,\gamma \in\Z_{+}^n,  \epsilon,\epsilon' \in\Z_2^n,
w\in W\}$ form a linear basis for~$\bH_W$ (the PBW basis).
\end{Theorem}
\begin{proof}
Recall that $W$ acts diagonally on $V = \h^{*} \oplus\h$. The
strategy of proving the theorem follows the suggestion of
\cite{W1} to modify \cite[Proof of Theorem~1.3]{EG} as follows.

Clearly $K := \Cl_{2n} \rtimes \C W$ is a semisimple algebra. Observe
that $E := V \otimes_{\C} K$ is a natural $K$-bimodule (even
though $V$ is not) with the right $K$-module structure on $E$
given by right multiplication and the left $K$-module structure on
$E$ by letting
\begin{gather*}
w \cdot (v\otimes a)  =  v^{w}\otimes wa,\\
c_i \cdot (x_j \otimes a)  =  (-1)^{\delta_{i j}} x_j \otimes (c_i
a), \\
c_i \cdot (y_j \otimes a)  =  y_j \otimes (c_i a),\\
e_i \cdot (x_j \otimes a)  =   x_j \otimes (e_i a), \\
e_i \cdot (y_j \otimes a)  =  (-1)^{\delta_{i j}} y_j \otimes (e_i a),
\end{gather*}
where  $v\in V$, $w\in W$, $a\in K$.

The rest of the proof can proceed in the same way as in
\cite[Proof of Theorem~1.3]{EG}, and it boils down to the
verif\/ications of the conjugation invariance (by $c_i$, $e_i$ and $W$) of
the def\/ining relations (\ref{Ayixj})--(\ref{Ayixi}),
(\ref{Dyixj})--(\ref{Dyixi}), or (\ref{Byixj})--(\ref{Byixi}) for type
$A$, $D$ or $B$ respectively, and the verif\/ication of the Jacobi
identities among the generators $x_i$ and $y_i$ for $1\le i \le n$.

Such verif\/ications are left to Lemmas~\ref{conj-inv-c},
\ref{conj-inv-W} and \ref{Jacobi} below.
\end{proof}

\begin{Remark}\rm
The algebra $\bH_W$ has two dif\/ferent triangular decompositions:
\begin{gather*}
\bH_W  \cong \C[\h^*] \otimes (\Cl_{2n} \rtimes\C W) \otimes \C[\h], \\
\bH_W  \cong (\C[\h^*] \otimes \Cl_{\h^*}) \otimes \C W \otimes
(\Cl_{\h} \otimes \C[\h]).
\end{gather*}
\end{Remark}

The detailed proofs of Lemmas~\ref{conj-inv-c}, \ref{conj-inv-W}
and \ref{Jacobi} below (also compare \cite{KW2}) are postponed to
the Appendix.

\begin{Lemma} \label{conj-inv-c}
Let $W=W_{A_{n-1}}, W_{D_n}$ or $W_{B_n}$. Then the relations
\eqref{Ayixj}--\eqref{Ayixi}, \eqref{Dyixj}--\eqref{Dyixi}, or
\eqref{Byixj}--\eqref{Byixi} are invariant under the conjugation by
$c_i$ and $e_i$ respectively, $1\le i \le n$.
\end{Lemma}

\begin{Lemma} \label{conj-inv-W}
The relations \eqref{Ayixj}--\eqref{Ayixi},
\eqref{Dyixj}--\eqref{Dyixi}, or \eqref{Byixj}--\eqref{Byixi}
are invariant under the conjugation by elements in
$W_{A_{n-1}}$, $W_{D_n}$ or $W_{B_n}$ respectively.
\end{Lemma}

\begin{Lemma} \label{Jacobi}
Let $W=W_{A_{n-1}}$, $W_{D_n}$ or $W_{B_n}$. Then the Jacobi
identity holds for any triple among $x_i$, $y_i$ in
$\bH_W$ for $1\le i \le n$.
\end{Lemma}

\begin{Remark} \label{automorphism}\rm
For $W = W_{A_{n-1}}$, $W_{D_n}$ or $W_{B_n}$, the algebra $\bH_W$
has a natural superalgebra structure by     letting $x_i$, $y_i$,
$s_j$ be even and $c_k$, $e_k$ be odd for all admissible $i$, $j$, $k$.
Moreover, the map $\varpi: \bH_W \longrightarrow \bH_W$ which
sends
\[
x_i \mapsto y_i, \qquad y_i \mapsto -x_i, \qquad c_i \mapsto e_i,
\qquad e_i \mapsto -c_i, \qquad s_j \mapsto s_j \qquad \forall\, i,j
\]
is an automorphism of $\bH_W$.
\end{Remark}

\subsection{The Dunkl representations}
Recall $K = \Cl_{2n} \rtimes \C W$. Denote by $\aH_y$ the
subalgebra of $\bH_W$ generated by $K$ and $y_1,\ldots,y_n$. A~$K$-module $M$ can be extended to $\aH_y$-module by demanding the
action of each $y_i$ to be trivial. We def\/ine
\[
M_y:=\text{Ind}_{\aH_y}^{\bH_W} M.
\]
Under the identif\/ication of vector spaces
\[
M_y = \C[x_1,\ldots,x_n]\otimes M,
\]
the action of $\bH_W$ on $M_y$ is transferred to
$\C[x_1,\ldots,x_n]\otimes M$ as follows. K acts on
$\C[x_1,\ldots,x_n]\otimes M$ by the following formulas:
\begin{gather*}
w \cdot (x_j\otimes m)  = x_j^{w}\otimes w m,\\
c_i \cdot (x_j \otimes m)  =  (-1)^{\delta_{i j}} x_j \otimes c_i m,\\
e_i \cdot (x_j \otimes m)  = x_j \otimes e_i m,
\end{gather*}
where $c_i,e_i \in \Cl_{2n}$, $w \in W$. Moreover, $x_i$ acts by
left multiplication in the f\/irst tensor factor, and the action of
$y_i$ will be given by the so-called Dunkl operators which we
compute below (compare~\cite{Dun, DO}).

A simple choice for a $K$-module is $\Cl_{2n}$, whose $K$-module
structure is def\/ined by letting $\Cl_{2n}$ act by left
multiplication and $W$ act diagonally.

\subsubsection[The Dunkl Operators for type $A_{n-1}$ case]{The Dunkl Operators for type $\boldsymbol{A_{n-1}}$ case}

We f\/irst prepare a few lemmas. It is understood in this paper that
the ratios of two (possibly noncommutative) operators $g$ and $h$
always means that $\frac{h}{g} =\frac1{g} \cdot h$.

\begin{Lemma} \label{Acomm:[y,x^l]}
Let $W=W_{A_{n-1}}$. Then the following holds in $\bH_W$ for $l \in \Z_+$ and $i\neq j$:
\begin{gather*}
\lbrack y_{i},x_{j}^l]  =  u\left(\frac{x_j^l - x_i^l}{x_j-x_i}
+ \frac{x_j^l -(-1)^l x_i^l}{x_j+x_i}c_ic_j\right)(1-e_ie_j)s_{ij},\\
\lbrack y_{i},x_{i}^l]  =  t c_i e_i \frac{x_i^l-(-x_i)^l}{2x_i}
-u
\sum_{k\neq i} \left( \frac{x_i^l - x_k^l}{x_i - x_k} + \frac{x_i^l
- (-x_k)^l}{x_i + x_k}c_k c_i\right)(1+e_k e_i)s_{k i}.
\end{gather*}
\end{Lemma}

\begin{proof}
This lemma is a type $A$ counterpart of Lemma~\ref{Bcomm:[y,x^l]}
for type $B$ below. A proof can be simply obtained by modifying
the proof of Lemma~\ref{Bcomm:[y,x^l]} with the removal of those
terms involving $\overline{s}_{ij}$, $\overline{s}_{ki}$, $\tau_i$
therein.
\end{proof}

\begin{Lemma} \label{Acomm:[y,f]}
Let $W=W_{A_{n-1}}$, and $f \in \C[x_1,\ldots,x_n]$. Then the
following identity holds in $\bH_W$:
\[
\lbrack y_{i},f] = tc_i e_i \frac{f-f^{\tau_i}}{2x_i}  -u
\sum_{k\neq i} \left( \frac{f -f^{s_{ki}}}{x_i - x_k} + \frac{f c_k
c_i - c_k c_i f^{s_{k i}}}{x_i+ x_k}\right)(1+e_k e_i)s_{k i}.
\]
\end{Lemma}

\begin{proof}
It suf\/f\/ices to check the formula for every monomial $f$ of the
form $x_1^{l_1}\cdots x_n^{l_n}$, which follows by
Lemma~\ref{Acomm:[y,x^l]} and an induction on $a$ based on the
identity
\begin{gather*}
 \lbrack y_{i},x_1^{l_1} \cdots x_a^{l_a}x_{a+1}^{l_{a+1}}]
 = \lbrack y_{i},x_1^{l_1} \cdots x_a^{l_a}]x_{a+1}^{l_{a+1}}
 + x_1^{l_1} \cdots x_a^{l_a} \lbrack y_{i}, x_{a+1}^{l_{a+1}}].
\tag*{\qed}
  \end{gather*}
  \renewcommand{\qed}{}
\end{proof}

Now we are ready to compute the Dunkl operator for $y_i$.

\begin{Theorem}\label{ADunkl:y}
Let $W=W_{A_{n-1}}$ and $M$ be a $(\Cl_{2n}\rtimes \C W)$-module.
The action of $y_i$ on the module $\C[x_1,\ldots,x_n] \otimes M$
is realized as the following Dunkl operators: for any $f
\in\C[x_1,\ldots,x_n]$ and $m \in M$, we have
\[
y_i \circ (f \otimes m) =  tc_i e_i \frac{f-f^{\tau_i}}{2x_i}
\otimes m - u \sum_{k\neq i} \left( \frac{f -f^{s_{ki}}}{x_i - x_k}
+ \frac{f c_k c_i - c_k c_i f^{s_{k i}}}{x_i+ x_k}\right)
\otimes(1+e_k e_i) s_{k i} m.
\]
\end{Theorem}
\begin{proof}
We calculate that
\[
y_i \circ (f \otimes m) = [y_i,f]\otimes m + f \otimes y_i m =
[y_i,f]\otimes m.
\]
Now the result follows from Lemma ~\ref{Acomm:[y,f]}.
\end{proof}

\subsubsection[The Dunkl Operators for type ${B_{n}}$ case]{The Dunkl Operators for type $\boldsymbol{{B_{n}}}$ case}

The proofs of Lemmas~\ref{Bcomm:[y,x^l]} and \ref{Bcomm:[y,f]} are
postponed to the Appendix.

\begin{Lemma} \label{Bcomm:[y,x^l]}
Let $W=W_{B_n}$. Then the following holds in $\bH_W$ for $l \in \Z_+$ and $i\neq
j$:
\begin{gather*}
\lbrack y_{i},x_{j}^l]  =  u\left(\frac{x_j^l - x_i^l}{x_j-x_i}
+ \frac{x_j^l -(-1)^l x_i^l}{x_j+x_i}c_ic_j\right)(1-e_ie_j)s_{ij}\\
\phantom{\lbrack y_{i},x_{j}^l]  =}{}  - u \left(\frac{x_j^l - (-x_i)^l}{x_j + x_i} - \frac{x_j^l
- x_i^l}{x_j - x_i}c_i c_j\right)(1+e_ie_j)\overline{s}_{i j}, \\
\lbrack y_{i},x_{i}^l]  =  t c_i e_i \frac{x_i^l-(-x_i)^l}{2x_i}
 -u
\sum_{k\neq i} \left( \frac{x_i^l -x_k^l}{x_i - x_k}
+ \frac{x_i^l - (-x_k)^l}{x_i + x_k}c_kc_i\right)(1+e_k e_i)s_{k i}\\
\phantom{\lbrack y_{i},x_{i}^l]  =}{}   - u \sum_{k\neq i} \left(\frac{x_i^l - (-x_k)^l}{x_i + x_k}
-\frac{x_i^l - x_k^l}{x_i - x_k}c_k
c_i\right)(1-e_ke_i)\overline{s}_{ki}
- v \frac{x_i^l -
(-x_i)^l}{2x_i}\tau_i.
\end{gather*}
\end{Lemma}

\begin{Lemma} \label{Bcomm:[y,f]}
Let $W=W_{B_n}$, and $f \in \C[x_1,\ldots,x_n]$. Then the
following holds in $\bH_W$:
\begin{gather*}
\lbrack y_{i},f]  =  tc_i e_i \frac{f-f^{\tau_i}}{2x_i}-u
\sum_{k\neq i} \left( \frac{f -f^{s_{ki}}}{x_i - x_k}
+ \frac{f -  f^{\overline{s}_{ki}}}{x_i+ x_k}c_k c_i\right)(1+e_k e_i)s_{k i} \\
\phantom{\lbrack y_{i},f]  =}{}  - v \frac{f -f^{\tau_i}}{2x_i} \tau_i - u \sum_{k\neq i}
\left(\frac{f - f^{\overline{s}_{k i}}}{x_i + x_k}
- \frac{f - f^{s_{k i}}}{x_i - x_k}c_k c_i\right)(1-e_ke_i)\overline{s}_{k i}.
\end{gather*}
\end{Lemma}
Now we are ready to compute the Dunkl operator for $y_i$.

\begin{Theorem}\label{BDunkl:y}
Let $W=W_{B_{n}}$ and $M$ be a $(\Cl_{2n}\rtimes \C W)$-module.
The action of $y_i$ on the module $\C[x_1,\ldots,x_n] \otimes M$
is realized as the following Dunkl operators: for any $f \in
\C[x_1,\ldots,x_n]$ and $m \in M$, we have
\begin{gather*}
y_i \circ (f \otimes m)  =   tc_i e_i \frac{f-f^{\tau_i}}{2x_i}
\otimes m
- u \sum_{k\neq i} \left( \frac{f -f^{s_{ki}}}{x_i - x_k}
+ \frac{f - f^{\overline{s}_{k i}}}{x_i+ x_k}c_k c_i\right)
\otimes(1+e_k e_i) s_{k i} m\\
\phantom{y_i \circ (f \otimes m)  =}{}  - u \sum_{k\neq i} \left(\frac{f - f^{\overline{s}_{k i}}}{x_i
+x_k} - \frac{f - f^{s_{k i}}}{x_i - x_k}c_k c_i\right) \otimes
(1-e_k e_i)\overline{s}_{k i}m
- v \frac{f - f^{\tau_i}}{2x_i}
\otimes \tau_i m.
\end{gather*}
\end{Theorem}
\begin{proof}
We observe that
\[
y_i \circ (f \otimes m) = [y_i,f]\otimes m + f \otimes y_i m =
[y_i,f]\otimes m.
\]
Now the result follows from Lemma ~\ref{Bcomm:[y,f]}.
\end{proof}

\subsubsection[The Dunkl Operators for type $D_{n}$ case]{The Dunkl Operators for type $\boldsymbol{D_{n}}$ case}

Due to the similarity of the bracket relations $\lbrack -,-]$ in
$D_{n}$ and $B_{n}$ cases (e.g.\ compare (\ref{Dyixi}) with
(\ref{Byixi})), the formula below for type $D_n$ is obtained from
its type $B_n$ counterpart in the previous subsection by dropping
the terms involving the parameter $v$.

\begin{Theorem} 
Let $W=W_{D_n}$, and let $M$ be a $(\Cl_{2n}\rtimes \C W)$-module.
The action of $y_i$ on $\C[x_1,\ldots,x_n] \otimes M$ is realized
as the following Dunkl operators: for any $f \in
\C[x_1,\ldots,x_n]$ and $m \in M$, we have
\begin{gather*}
y_i \circ (f \otimes m)  = tc_i e_i \frac{f-f^{\tau_i}}{2x_i}
\otimes m
- u \sum_{k\neq i} \left( \frac{f -f^{s_{ki}}}{x_i - x_k}
+ \frac{f  - f^{\overline{s}_{k i}}}{x_i+ x_k}c_k c_i\right)
\otimes(1+e_k e_i) s_{k i} m\\
\phantom{y_i \circ (f \otimes m)  =}{}  - u \sum_{k\neq i} \left(\frac{f - f^{\overline{s}_{k i}}}{x_i
+ x_k} - \frac{f - f^{s_{k i}}}{x_i - x_k}c_k c_i\right)\otimes
(1-e_k e_i)\overline{s}_{k i}m.
\end{gather*}
\end{Theorem}

\subsection[The even center for $\bH_W$]{The even center for $\boldsymbol{\bH_W}$}

Recall that the {\em even center} $\mathcal Z (A)$ of a
superalgebra $A$ consists of the even central elements of~$A$. It
turns out the algebra $\bH_W$ has a large center.

\begin{Proposition} \label{CenDaHa}
Let $W$ be either $W_{A_{n-1}}$, $W_{D_n}$ or $W_{B_n}$. The even
center $\Zn(\bH_W)$ contains $\C[x_1^2,\ldots, x_n^2]^W$ and $\C
[y_1^2,\ldots, y_n^2]^W$ as subalgebras. In particular, $\bH_W$ is
module-finite over its even center.
\end{Proposition}
\begin{proof}
Let $f\in\C[x_1^2,\ldots, x_n^2]^W$. Then $f-f^{\tau_i}=0$ for each
$i$. Moreover, by the def\/inition of~$\bH_W$, $f$~commutes with
$\Cl_{2n}$, $W$, and $x_i$ for all $1\leq i \leq n$. Since $f =
f^{w}$ for all $w \in W$, it follows from Lemmas \ref{Acomm:[y,f]}
and \ref{Bcomm:[y,f]} that $\lbrack y_i, f] = 0$ for each $i$. Hence
$f$ commutes with $\C[y_1,\ldots,y_n]$. Therefore $f$ is in the even
center $\Zn (\bH_W)$. It follows from the automorphism $\varpi$ of
$\bH_W$ def\/ined in Remark~\ref{automorphism} that $\C[y_1^2,\ldots,
y_n^2]^W$ must also be in the even center $\Zn (\bH_W)$.
\end{proof}

%
%
%
%
\section[The spin double affine Hecke-Clifford algebras]{The spin double af\/f\/ine Hecke--Clif\/ford algebras} \label{sec:sDaHCa}

Recall that $W$ is one of the classical Weyl groups of type
$A_{n-1}$, $B_n$, or $D_n$. The goal of this section is to introduce
and study the spin double af\/f\/ine Hecke--Clif\/ford algebra (sDaHCa)
$\sH_W$, which is, roughly speaking, obtained by decoupling the
Clif\/ford algebra $\Cl_\h$ from the DaHa $\bH_W$ in
Section~\ref{sec:DaHa:doubleClifford}. The spin Weyl group algebra
$\C W^-$ appears naturally in the process. We remark that the
algebra $\sH_W$ is dif\/ferent from either the spin double af\/f\/ine
Hecke algebra or the double af\/f\/ine Hecke--Clif\/ford algebra introduced
in \cite{W1, KW2}.

\subsection[The definition of sDaHCa $\sH_W$]{The def\/inition of sDaHCa $\boldsymbol{\sH_W}$}

Following \cite{KW2}, we introduce the notation
\begin{gather*}
t_{i\uparrow j}  =  \left \{
 \begin{array}{ll}
 t_it_{i+1}\cdots t_j, & \text{ if } i\leq j,\\
 1, &\text{ otherwise},
 \end{array}
 \right. \qquad
t_{i\downarrow j}  =  \left \{
 \begin{array}{ll}
 t_it_{i-1}\cdots t_j, & \text{ if } i\geq j,\\
 1, &\text{ otherwise}.
 \end{array}
 \right.
\end{gather*}
Def\/ine the following odd elements in $\C W^-$ of order $2$, which
are an analogue of ref\/lections in~$W$, for $1\leq i <j\leq n$:
\begin{gather*}
 t_{i j} \equiv [i,j]
 = (-1)^{j-i-1}t_{j-1}\cdots t_{i+1}t_{i}t_{i+1}\cdots t_{j-1},\\
 t_{ji} \equiv {[}j,i] = -[i,j], \\
\overline{t}_{i j} \equiv \overline{[i,j]}
 = \left \{
 \begin{array}{ll}
  (-1)^{j-i-1} t_{j \uparrow n-1} t_{i \uparrow n-2}
  {t}_{n}
   t_{n-2 \downarrow i} t_{n-1 \downarrow j}, & \text{for type } D_n, \\
  (-1)^{j-i} t_{j \uparrow n-1} t_{i \uparrow n-2}
   {t}_{n}{t}_{n-1}{t}_n
   t_{n-2 \downarrow i} t_{n-1 \downarrow j}, & \text{for type } B_n,
  \end{array}
  \right.
\\
\overline{t}_{ji} \equiv \overline{[j,i]}  =  \overline{[i,j]}, \\
 \overline{t}_i \equiv \overline{[i]}
  =  (-1)^{n-i} t_i\cdots t_{n-1}t_n t_{n-1}\cdots t_i \qquad (1\le i \leq n).
\end{gather*}
 The notations $[i,j]$, $\overline{[i,j]}$ here are consistent with
the inclusions of algebras $\C W_{A_{n-1}}^- \leq \C W_{D_n}^-
\leq \C W^-_{B_n}$.
%
%

As in \cite{W1} (also cf.~\cite{KW1, KW2}), a {\em skew-polynomial
algebra} is the $\C$-algebra generated by $b_1,\ldots,b_n$ subject
to the relations $ b_ib_j+b_jb_i =0,\quad (i\neq j).$ This
algebra, denoted by $\Cl[b_1,\ldots,b_n]$, is naturally a
superalgebra by letting each $b_i$ be odd, and it has a linear
basis given by $b^\alpha :=b_1^{k_1}\cdots b_n^{k_n}$ for $\alpha
=(k_1,\ldots,k_n) \in \Z_+^n$.

Consider the group homomorphism $\rho: W_{B_n} \rightarrow S_n$
def\/ined by $\rho(s_i)=s_i$ and $\rho(s_n) = 1$ for $1\le i \le
n-1$. By restriction if needed, we have a group homomorphism
\[
\rho: W \longrightarrow S_n, \qquad \sigma \mapsto \rho(\sigma) =
\sigma^*
\]
for $W=W_{A_{n-1}}, W_{B_n}$ or $ W_{D_n}$. Observe that $\tau_i^*
= 1$ and $\overline{s}_{ij}^* = s_{ij}$ for all $1\le i\neq j \le
n$.

\begin{Definition}\rm
Let $t,u,v \in \C$, and $W =W_{A_{n-1}}, W_{D_{n}}$, or
$W_{B_{n}}$. The sDaHCa $\sH_W$ is the algebra generated by
$x_i$, $\eta_i$ $(1\leq i \leq n)$ and $\Cl_{\h^*} \rtimes \C W^-$,
subject to the relations
\begin{gather*}
\eta_i \eta_j  = -\eta_j \eta_i, \quad x_i x_j = x_j x_i  \qquad (i \neq j),\\
c_i \eta_j  = - \eta_j c_i, \quad c_i x_j
= (-1)^{\delta_{ij}}x_j c_i  \qquad(\forall\, i,j),\\
t_i x_j  = x_j^{s_i}t_i, \quad t_i \eta_j = -\eta_j^{s_i^*}t_i
\qquad (t_i \in \C W^-)
\end{gather*}
and the following additional relations:
\begin{gather*}
\text{Type A:}
 \ \
\begin{cases}
\lbrack \eta_i, x_j] = u (1+c_ic_j)[i,j] \qquad (i\neq j), \\
\displaystyle \lbrack \eta_i, x_i] = t c_i + u \sum_{k\neq i} (1+c_kc_i)[k,i],
\end{cases}\\
\text{Type D:}
 \ \
\begin{cases}
\lbrack \eta_i, x_j]
= u((1+c_i c_j)[i,j]-(1-c_i c_j)\overline{[i,j]}) \qquad (i\neq j), \\
\displaystyle \lbrack \eta_i, x_i] = t  c_i +u \sum_{k\neq i} \big((1+c_k
c_i)[k,i]-(1-c_k c_i)\overline{[k,i]}\big),
\end{cases}\\
\text{Type B:} \ \
\begin{cases}
\lbrack \eta_i, x_j]
= u((1+c_i c_j)[i,j]-(1-c_i c_j)\overline{[i,j]}) \qquad (i\neq j), \\
\displaystyle\lbrack \eta_i, x_i] = t  c_i +u \sum_{k\neq i} \big((1+c_k
c_i)[k,i]-(1-c_k c_i)\overline{[k,i]}\big)+ v\overline{[i]}.
\end{cases}
\end{gather*}
\end{Definition}

\subsection{Isomorphism of superalgebras}

For $W=W_{A_{n-1}}, W_{B_n},$ or $W_{D_n}$, we recall an algebra
isomorphism (see \cite[Lemma 5.4]{KW2})
\[
    \Phi: \  \Cl_{\h} \rtimes \C W \rightarrow \Cl_{\h} \otimes \C W^-
\]
which sends
\begin{gather}
(e_k-e_i) s_{ik}  \longmapsto -\sqrt{-2}\;[k,i], \nonumber \\
(e_k +e_i) \overline{s}_{ik}  \longmapsto
-\sqrt{-2}\;\overline{[k,i]},  \label{transpD}\\
e_i \tau_i  \longmapsto -\sqrt{-1} \, \overline{[i]}\nonumber
\end{gather}
for $i \neq k$, whenever it is applicable. The inverse of $\Phi$
is denoted by $\Psi$.

Note that the algebra $\sH_W$ has a natural superalgebra structure
by letting each $\eta_i$, $c_i$, $t_j$ be odd and $x_i$ be even for
all admissible $i$, $j$.
\begin{Theorem} \label{th:isomADBsdahc}
Let $W$ be $W_{A_{n-1}}$, $W_{D_n}$ or
$W_{B_n}$. Then,
\begin{enumerate}\itemsep=0pt
\item[$1)$] there exists an isomorphism of superalgebras
\[
\Phi: \  \bH_W(t,u,v)\longrightarrow\Cl_{\h}\otimes \sH_W(-t,-\sqrt{-2}u,\sqrt{-1}v)
\]
which extends $\Phi: \Cl_{\h} \rtimes \C W \rightarrow \Cl_{\h} \otimes
\C W^-$ and sends
\[
y_i \mapsto e_i \eta_i,\qquad
x_i \mapsto x_i,\qquad
c_i \mapsto c_i,  \qquad \forall\, i;
\]

\item[$2)$] the inverse
\[
\Psi: \ \Cl_{\h}\otimes \sH_W(-t,-\sqrt{-2}u,\sqrt{-1}v) \longrightarrow \bH_W(t,u,v)
\]
extends $\Psi:\Cl_{\h} \otimes\C W^- \rightarrow \Cl_{\h} \rtimes \C W$ and sends
\[
\eta_i \mapsto e_i y_i,\qquad
x_i \mapsto x_i,\qquad
c_i \mapsto c_i,  \qquad \forall \, i.
\]
\end{enumerate}
\end{Theorem}
\begin{proof}
We need to check that $\Phi$ preserves the relations
(\ref{commonA}), (\ref{Ayixj})--(\ref{Ayixi}),
(\ref{Dyixj})--(\ref{Dyixi}), and (\ref{Byixj})--(\ref{Byixi}) for
$W=W_{A_{n-1}}, W_{D_n}$, and $W_{B_n}$ respectively.

First, we shall verify that $\Phi$ preserves
(\ref{Byixj})--(\ref{Byixi}) with $W=W_{B_n}$. Indeed, by
\eqref{transpD} or \cite[Lemma 5.4]{KW2}, we have
\begin{gather*}
\Phi (\text{l.h.s. of (\ref{Byixj})}) = e_i [\eta_i, x_j] =-\sqrt{-2}u e_i \big((1-c_j c_i)[i,j]-(1+c_j c_i)\overline{[i,j]}\big)\\
\qquad{} = \Phi \big(u\big( (1+c_{i}c_{j})(1+e_{j}e_{i})s_{ji}-(1-c_{i}
c_{j})(1-e_{j}e_{i})\overline{s}_{ij}\big ) \big ) \\
\qquad = \Phi (\text{r.h.s. of (\ref{Byixj})}).
\end{gather*}

Also, we have
\begin{gather*}
 \Phi (\text{l.h.s. of (\ref{Byixi})}) = e_i [\eta_i, x_i] \\
\qquad {}=-t \cdot e_i c_i - \sqrt{-2}u e_i\sum_{k\neq i}
\big((1+c_k c_i)[k,i]-(1-c_k c_i)\overline{[k,i]}\big)+ \sqrt{-1}ve_i\overline{[i]} \\
\qquad {} =\Phi\big(tc_i e_i -u\sum_{k\neq i}
\big((1+c_{k}c_{i})(1+e_k e_i)s_{ki}+(1-c_{k} c_{i})(1-e_ke_i)\overline{s}_{ki}\big)
-v \tau_i\big)\\
\qquad{}= \Phi (\text{r.h.s. of (\ref{Byixi})}).
\end{gather*}
It is easy to check that $\Phi$ preserves (\ref{commonA}), and we
will restrict ourselves to verify just a few relations among
(\ref{commonA}). For $j\neq i, i+1$, we have
\[
    \Phi(s_i y_j) = -\sqrt{-1}\nu_it_i e_j \eta_j
    = -\sqrt{-1}e_j\eta_j\nu_i t_i = \Phi(y_j s_i).
\]
Moreover,
\[
\Phi(s_n y_n)  = -\sqrt{-1}\nu_n t_n e_n \eta_n = \sqrt{-1}t_n
\eta_n = -\sqrt{-1}\eta_n t_n = \Phi(-y_n s_n).
\]
This proves that $\Phi$ is an algebra homomorphism for type $B_n$.

By dropping the terms involving $v$ in the above equations, we
verify that the relations (\ref{Dyixj})--(\ref{Dyixi}) with
$W=W_{D_n}$ are preserved by $\Phi$. By further dropping the terms
involving $\overline{[ij]}$, $\overline{s}_{ij}$ etc., we also
verify (\ref{Ayixj})--(\ref{Ayixi}) with $W =W_{A_{n-1}}$. So, the
homomorphism $\Phi$ is well def\/ined in all cases.

Similarly, one shows that $\Psi$ is a well-def\/ined algebra
homomorphism. Since $\Phi$ and $\Psi$ are inverses on generators,
they are (inverse) algebra isomorphisms.
\end{proof}

The isomorphism in Theorem \ref{th:isomADBsdahc} exactly means
that the superalgebras $\bH_W$ and $\sH_W$ are Morita
super-equivalent in the sense of \cite{W2}.

\begin{Corollary} \label{CensDaHCa}
Let $W$ be one of the Weyl groups $W_{A_{n-1}}$, $W_{D_n}$ or
$W_{B_n}$. The even center $\Zn(\sH_W)$ of $\sH_W$ contains
$\C[\eta_1^2,\ldots,\eta_n^2]^W$ and $\C [x_1^2,\ldots,x_n^2]^W$.
In particular, $\sH_W$ is module-finite over its even center.
\end{Corollary}
\begin{proof}
By the isomorphism $\Phi$ in Theorem \ref{th:isomADBsdahc} and the
Proposition \ref{CenDaHa}, we have that $\Zn(\Cl_{\h} \otimes
\sH_W)$ contains the subalgebras $\C [\eta_1^2,\ldots,\eta_n^2]^W$
and $\C [x_1^2,\ldots,x_n^2]^W$, and so does $\Zn(\sH_W)$.
\end{proof}

\subsection[The PBW property for $\sH_W$]{The PBW property for $\boldsymbol{\sH_W}$}

We have the following PBW type property for the algebra $\sH_W$.
\begin{Theorem} \label{PBW:sdahc}
Let $W$ be one of the Weyl groups $W_{A_{n-1}}$, $W_{D_n}$ or
$W_{B_n}$. The multiplication of the subalgebras induces an
isomorphism of vector spaces
\[
\Cl[\eta_1,\ldots,\eta_n] \otimes\Cl_{\h^*} \otimes \C W^- \otimes\C
[\h^*] \longrightarrow \sH_W.
\]
Equivalently, the set $\{\eta^{\al} c^{\epsilon}\sigma
x^{\gamma}\}$ forms a basis for $\sH_W$, where $\sigma$ runs over
a basis for $\C W^-$, $\epsilon \in \Z_2^n$, and $\al, \gamma \in
\Z_+^n$.
\end{Theorem}

\begin{proof}
It follows from the def\/ining relations that $\sH_W$ is spanned by
the elements $\eta^{\al} c^{\epsilon} \sigma x^{\gamma}$ where
$\sigma$ runs over a basis for $\C W^-$, $\al, \gamma \in \Z_+^n$,
and $\epsilon \in \Z_2^n$. By the isomorphism $\Psi:
\Cl_{\h}\otimes \sH_W\longrightarrow\bH_W$ in
Theorem~\ref{th:isomADBsdahc}, we see that the image $\Psi
(\eta^{\al} c^{\epsilon}\sigma x^{\gamma})$ are linearly
independent in $\bH_W$ by the PBW property for $\bH_W$ (see
Theorem~\ref{PBW:ADBdahc}). So the elements $\eta^{\al}
c^{\epsilon}\sigma x^{\gamma}$ are linearly independent in
$\sH_W$. Therefore, the set $\{\eta^{\al} c^{\epsilon} \sigma
x^{\gamma}\}$ forms a basis for $\sH_W$.
\end{proof}

\subsection[The Dunkl operators for $\sH_W$]{The Dunkl operators for $\boldsymbol{\sH_W}$}

Denote by $\h_{\eta}$ the subalgebra of $\sH_W$ generated by
$\eta_i$ $(1\leq i \leq n)$ and $\Cl_{\h^*} \rtimes_- \C W^-$. A
$(\Cl_{\h^*} \rtimes_- \C W^-)$-module $V$ can be extended to a
$\h_{\eta}$-modules by letting the actions of $\eta_i$ on $V$ to
be trivial for each $i$. We def\/ine
\[
    V_{\eta}:=\text{Ind}_{\h_{\eta}}^{\sH_W} V \cong \C[x_1,\ldots,x_n]\otimes V.
\]
On $\C[x_1,\ldots,x_n]\otimes V$, the element $t_i \in \C W^-$
acts as $s_i \otimes t_i$, $c_i \in \Cl_{\h^*}$ acts by $c_i \cdot
(x_j \otimes v) = (-1)^{\delta_{ij}}x_j \otimes c_i v$, and $x_i$
acts by left multiplication, and $\eta_i$ acts as anti-commuting
Dunkl operators, which we will describe in this section.

Under the superalgebra isomorphism
$\Phi:\bH_W\rightarrow\Cl_n\otimes \sH_W $ in
Theorem~\ref{th:isomADBsdahc}, we obtain anti-commuting Dunkl
operators $\eta_i$ by fairly straightforward computation. They are
counterparts of those in Section~\ref{sec:DaHa:doubleClifford},
and we omit the proofs.

\subsubsection[Dunkl operator for type $A_{n-1}$]{Dunkl operator for type $\boldsymbol{A_{n-1}}$}

The following is a counterpart of Theorem~\ref{ADunkl:y}.

\begin{Proposition}
Let $W=W_{A_{n-1}}$ and $V$ be a $\Cl_n\rtimes \C W^-$-module. The
action of $\eta_i$ on the $\sH_W$-module $\C[x_1,\ldots,x_n]
\otimes V$ is realized as a Dunkl operator as follows. For any
polynomial $f \in \C[x_1,\ldots,x_n]$ and $m \in V$, we have
\begin{gather*}
\eta_i \circ (f \otimes m)  =   tc_i \frac{f-f^{\tau_i}}{2x_i}
\otimes m
+ u \sum_{k\neq i} \left( \frac{f -f^{s_{ki}}}{x_i - x_k}
+ \frac{f c_k c_i - c_k c_k f^{s_{k i}}}{x_i+ x_k}\right)
\otimes[k,i] m.
\end{gather*}
\end{Proposition}

\subsubsection[Dunkl operator for type $B_{n}$]{Dunkl operator for type $\boldsymbol{B_{n}}$}

The following is a counterpart of Theorem~\ref{BDunkl:y}

\begin{Proposition} \label{BDunkl:etax}
Let $W=W_{B_{n}}$ and $V$ be a $(\Cl_n\rtimes \C W^-)$-module. The
action of $\eta_i$ on the $\sH_W$-module $\C[x_1,\ldots,x_n]
\otimes V$ is realized as a Dunkl operator as follows. For any
polynomial $f \in \C[x_1,\ldots,x_n]$ and $m \in V$, we have
\begin{gather*}
\eta_i \circ (f \otimes m)  =   tc_i \frac{f-f^{\tau_i}}{2x_i}
\otimes m
+ u \sum_{k\neq i} \left( \frac{f -f^{s_{ki}}}{x_i - x_k}
+ \frac{f - f^{\overline{s}_{k i}}}{x_i+ x_k}c_k c_i\right) \otimes[k,i] m\\
\phantom{\eta_i \circ (f \otimes m)  =}{}  - u \sum_{k\neq i} \left(\frac{f - f^{\overline{s}_{k i}}}{x_i +
x_k} - \frac{f - f^{s_{k i}}}{x_i - x_k}c_k c_i\right)\otimes
\overline{[k,i]}m
+ v \frac{f - f^{\tau_i}}{2x_i} \otimes
\overline{[i]} m.
\end{gather*}
\end{Proposition}

\subsubsection[Dunkl operator for type $D_{n}$]{Dunkl operator for type $\boldsymbol{D_{n}}$}

\begin{Proposition}
Let $W=W_{D_{n}}$ and $V$ be a $\Cl_n\rtimes \C W^-$-module. The
action of $\eta_i$ on the $\sH_W$-module $\C[x_1,\ldots,x_n]
\otimes V$ is realized as a Dunkl operator as follows. For any
polynomial $f \in \C[x_1,\ldots,x_n]$ and $m \in V$, we have
\begin{gather*}
\eta_i \circ (f \otimes m)  =   tc_i \frac{f-f^{\tau_i}}{2x_i}
\otimes m
+ u \sum_{k\neq i} \left( \frac{f -f^{s_{ki}}}{x_i - x_k}
+ \frac{f - f^{\overline{s}_{k i}}}{x_i+ x_k}c_k c_i\right) \otimes[k,i] m\\
\phantom{\eta_i \circ (f \otimes m)  =}{}  - u \sum_{k\neq i} \left(\frac{f - f^{\overline{s}_{k i}}}{x_i +
    x_k} - \frac{f - f^{s_{k i}}}{x_i - x_k}c_k c_i\right)\otimes
    \overline{[k,i]}m.
\end{gather*}
\end{Proposition}

\begin{Remark}\rm
The general formula for the Dunkl operators $\eta_i$ for $\sH_W$
resembles the Dunkl opera\-tor $y_i$ for $\dahc_W$ which appeared in
\cite[Theorems 4.4, 4.10, 4.14]{KW2}. However, $\eta_i \eta_j = -
\eta_j \eta_i$, while $y_i y_j = y_j y_i$ for $i\neq j$.
\end{Remark}

\section[The odd double affine Hecke algebras]{The odd double af\/f\/ine Hecke algebras}\label{sec:oddDaHa}

In this section, we shall introduce an  odd double af\/f\/ine Hecke
algebra $\HO_W$ which is generated by $\C W$ and two isomorphic
skew-polynomial subalgebras. Recall that $W$ is assumed to be one
of the classical Weyl groups of type $A_{n-1}$, $B_n$, or $D_n$.

Recall also the group homomorphism $\rho:W \longrightarrow S_n$
def\/ined in Section~\ref{sec:sDaHCa} which sends $\sigma \mapsto
\sigma^*$ for all $\sigma \in W$. We shall need two (isomorphic)
skew-polynomial algebras $\Cl \{\h^*\}=\Cl [\xi_1,\ldots,\xi_n]$
and  $\Cl \{\h\} =\Cl [\eta_1,\ldots,\eta_n]$, which are naturally
acted upon by the symmetric group $S_n$ or the group $W_{B_n}$
by permuting the indices possibly coupled with sign changes. We
shall denote the action of $\sigma\in W_{B_n}$ by $f \mapsto
f^\sigma$.

\subsection[The definition of $\HO_W$]{The def\/inition of $\boldsymbol{\HO_W}$}

As usual we denote $[\xi, \eta]_+ = \xi \eta + \eta \xi.$
\begin{Definition}\rm
Let $t,u,v \in \C$ and $W$ be $W_{A_{n-1}}, W_{D_{n}}$, or
$W_{B_{n}}$. The odd DaHa $\HO_W$ is the algebra generated by
$\xi_i$, $\eta_i$ $(1\leq i \leq n)$ and $\C W$, subject to the
relations
\begin{gather*}
\eta_i \eta_j  = -\eta_j \eta_i, \qquad \xi_i \xi_j = -\xi_j \xi_i  \qquad (i \neq j), \\
\sigma \xi_j = \xi_j^{\sigma^*}\sigma, \qquad \sigma \eta_j =
\eta_j^{\sigma^*} \sigma  \qquad (\sigma \in W)
\end{gather*}
and the following additional relations:
\begin{gather*}
 \text{Type A:} \ \
\begin{cases}
\lbrack \eta_i, \xi_j]_+ = u s_{ij} \qquad (i \neq j),\\
\displaystyle \lbrack \eta_i, \xi_i]_+ = t\cdot 1 + u\sum_{k \ne i} s_{ki},
\end{cases}\\
 \text{Type D:} \ \
\begin{cases}
\lbrack \eta_i, \xi_j]_+ = u \left(s_{ij}+\overline{s}_{ij}\right) \qquad (i \neq j),\\
\displaystyle \lbrack \eta_i, \xi_i]_+ = t\cdot 1 + u\sum_{k \ne i}\left(s_{ki} + \overline{s}_{ij}\right),
\end{cases}\\
 \text{Type B:} \ \
\begin{cases}
\lbrack \eta_i, \xi_j]_+
= u \left(s_{ij}+\overline{s}_{ij}\right) \qquad (i \neq j),\\
\displaystyle \lbrack \eta_i, \xi_i]_+
= t\cdot 1 + u\sum_{k \ne i}\left(s_{ki} + \overline{s}_{ij}\right)+ v \tau_i.
\end{cases}
\end{gather*}
\end{Definition}

The algebra $\HO_W$ has a natural superalgebra structure by
letting $s_j$ be even and $\eta_i$, $\xi_i$ be odd for all $i$, $j$.

\begin{Remark}\rm
The def\/ining relations for the algebra $\HO_W$ dif\/fer from those
for the usual rational DaHa (also known as rational Cherednik
algebra) $\daha_W$ \cite{EG} by signs. One can introduce a
so-called ``covering algebra'' $\wtd{\mathbb{H}}$ (as done in
\cite{W2, KW2} in similar setups) which contains a central element
$z$ of order $2$, so that the algebras $\daha_W$ and $\HO_W$ are
simply the quotients of $\wtd{\mathbb{H}}$ by the ideal generated
by $z-1$ and $z+1$ respectively.

The def\/inition of $\HO_W$ is motivated by the Morita
(super)equivalence with $\bH_W$ and $\sH_W$. The def\/ining
relations above suggest a further extension of odd DaHa associated
to the inf\/inite series complex ref\/lection groups.
\end{Remark}

\subsection{Isomorphism of superalgebras}

\begin{Lemma} \label{identity}
Let $W$ be one of the Weyl groups $W_{A_{n-1}}$, $W_{D_n}$ or
$W_{B_n}$. The isomorphism $\dot{\Phi}: \Cl_n \rtimes_- \C W^-
\rightarrow \Cl_n \otimes \C W$ (see
Theorem~{\rm \ref{th:isofinite:Cl_n}}) sends
\begin{gather*}
(c_k-c_i) [k,i]  \longmapsto \sqrt{2}\;s_{ki}, \qquad
%
(c_k +c_i)\overline{[k,i]}  \longmapsto
\sqrt{2}\; \overline{s}_{ik}, \qquad
%
c_i  \overline{[i]} \longmapsto  \tau_i.
\end{gather*}
\end{Lemma}
\begin{proof}
The lemma can be proved by induction  very similar to \cite[Lemma
5.4]{KW2}, and we skip the detail.
\end{proof}

\begin{Theorem}\label{th:isomADBodddaha}
Let $W$ be one of the Weyl groups $W_{A_{n-1}}$, $W_{D_n}$ or
$W_{B_n}$. Then,
\begin{enumerate}\itemsep=0pt
\item[$1)$] there exists an isomorphism of superalgebras
\[
\dot{\Phi}: \  \sH_W(t,u,v)\longrightarrow\Cl_n\otimes \HO_W(-t,\sqrt{2}u,-v)
\]
which extends $\dot{\Phi}: \Cl_n \rtimes_- \C W^- \rightarrow
\Cl_n \otimes \C W$ and sends
\[
\eta_i \mapsto \eta_i,\qquad
x_i \mapsto c_i\xi_i,  \qquad \forall\, i;
\]

\item[$2)$] the inverse
\[
\dot{\Psi}: \  \Cl_n\otimes \HO_W(-t,\sqrt{2}u,-v) \longrightarrow \sH_W(t,u,v)
\]
extends $\dot{\Psi}:\Cl_n \otimes\C W \rightarrow \Cl_n \rtimes_-
\C W^-$ and sends
\[
\eta_i \mapsto \eta_i,\qquad
\xi_i \mapsto c_i x_i,\qquad \forall\, i.
\]
\end{enumerate}
\end{Theorem}

\begin{proof}
We f\/irst need to check that $\dot{\Phi}$ preserves the def\/ining
relations of $\sH_W(t,u,v)$ and so $\Phi$ is a well-def\/ined
homomorphism. Using Lemma~\ref{identity}, we shall check a few
cases in type $B_{n}$ case, and leave the rest for the reader to
verify. For $i \neq j$, we have
    \begin{gather*}
\dot{\Phi}([\eta_i,x_j])  =-c_j [\eta_i,\xi_j]_+
 = -\sqrt{2} u c_j(s_{ij}+\overline{s}_{ij})\\
\phantom{\dot{\Phi}([\eta_i,x_j])}{} = \frac{u}{\sqrt{2}}\big((1+c_i c_j)(c_i - c_j)s_{ij}
 - (1-c_i c_j)(c_i + c_j)\overline{s}_{ij}  \big)\\
\phantom{\dot{\Phi}([\eta_i,x_j])}{}
=\dot{\Phi}\big(u\big((1+c_i c_j)[i,j]
- (1-c_i c_j)\overline{[i,j]}\big)\big),\\
\dot{\Phi}([\eta_i,x_i])  = -c_i [\eta_i,\xi_i]_+
= -c_i \Bigg(-t + \sqrt{2}u \sum_{k\neq i}\big(s_{ik}
+\overline{s}_{ik}\big) -v\tau_i \Bigg)\\
\phantom{\dot{\Phi}([\eta_i,x_i])}{} =t c_i  - \sqrt{2}u c_i \sum_{k\neq i}\big(s_{ik}
+\overline{s}_{ik}\big) +v c_i\tau_i\\
\phantom{\dot{\Phi}([\eta_i,x_i])}{}=\dot{\Phi}\Bigg(t c_i + u \sum_{k\neq i}\big((1+c_k
c_i)[k,i]+(1-c_k c_i)\overline{[k,i]}\big) + v \tau_i\Bigg).
    \end{gather*}
Also, if $j \neq n$, we have
    \begin{gather*}
\dot{\Phi}(t_n x_j)  = c_n s_n c_j \xi_j = c_j \xi_j c_n s_n = \dot{\Phi}(x_j t_n),\\
\dot{\Phi}(t_n x_n) = c_n s_n c_n \xi_n = s_n \xi_n = \dot{\Phi}(-x_n t_n),\\
\dot{\Phi}(t_n \eta_j)  = c_n s_n \eta_j = -\eta_j c_n s_n = \dot{\Phi}(-\eta_j t_n),\\
\dot{\Phi}(t_n \eta_n) = c_n s_n \eta_n = -\eta_n c_n s_n =
\dot{\Phi}(-\eta_n t_n).
    \end{gather*}

Similarly, one shows that $\dot{\Psi}$ is a well-def\/ined algebra
homomorphism. Since $\dot{\Phi}$ and $\dot{\Psi}$ are inverses on
generators, they are (inverse) algebra isomorphisms.
\end{proof}

The next corollary can be proved similarly to Corollary~\ref{CensDaHCa}.

\begin{Corollary}\label{CenOddDaHa}
Let $W$ be one of the Weyl groups $W_{A_{n-1}}$, $W_{D_n}$ or
$W_{B_n}$. The even center for $\HO_W$ contains $\C[\eta_1^2,\ldots,
\eta_n^2]^W$ and $\C [\xi_1^2,\ldots,\xi_n^2]^W$. In particular,
$\HO_W$ is module-finite over its even center.
\end{Corollary}

\begin{Example}\rm
Usually there are other central elements beyond those given in the
above corollary. For example, $\xi_1^2 \eta_2^2 + \xi_2^2 \eta_1^2
- u s_1 (\xi_1 -\xi_2)(\eta_1 -\eta_2)$ lies in $\Zn
(\HO_{W_{A_1}})$.
\end{Example}

\subsection[The PBW property for $\HO_W$]{The PBW property for $\boldsymbol{\HO_W}$}

We have the following PBW type property for the algebra $\HO_W$
which can be proved similarly to Theorem~\ref{PBW:sdahc}, using
now the isomorphism $\dot{\Phi}$.

\begin{Theorem}\label{PBW:ADBoddDaHa}
Let $W$ be one of the Weyl groups $W_{A_{n-1}}$, $W_{D_n}$ or
$W_{B_n}$. The multiplication of the subalgebras induces an
isomorphism of vector spaces
\[
\Cl[\xi_1,\ldots,\xi_n] \otimes \C W \otimes\Cl
[\eta_1,\ldots,\eta_n] \longrightarrow \HO_W.
\]
Equivalently, the set $\{\xi^{\al} \sigma \eta^{\gamma}\}$ forms a
basis for $\HO_W$, where $\sigma \in W$, and $\al, \gamma \in
\Z_+^n$.
\end{Theorem}

\subsection[The Dunkl operators for $\HO_W$]{The Dunkl operators for $\boldsymbol{\HO_W}$}
\label{subsec:Dunkl}

Denote by $\h_{\eta}$ the subalgebra of $\HO_W$ generated by
$\eta_i$ $(1\leq i \leq n)$ and $\C W$. Let $V$ be the trivial $\C
W$-module, and extend $V$ a $\h_{\eta}$-module by letting the
actions of each $\eta_i$ on V be trivial. Def\/ine
\[
    V_{\eta}:=\text{Ind}_{\h_{\eta}}^{\HO_W} V \cong \Cl[\xi_1,\ldots,\xi_n].
\]
On $\Cl[\xi_1,\ldots,\xi_n]$, $\sigma \in W$ acts as $\rho(\sigma)
=\sigma^*$, $\xi_i$ acts by left multiplication, and $\eta_i$ acts
as anti-commuting Dunkl operators which we establish below. (It is
easy to replace the trivial module above by any $\C W$-module.)

\subsubsection[Dunkl operator for type $A$ case]{Dunkl operator for type $\boldsymbol{A}$ case}

For each $i$, we introduce a super derivation $\partial_{\xi_i}$
on $\Cl[\xi_1,\ldots,\xi_n]$ def\/ined inductively by
$\partial_{\xi_i}(\xi_j) = \delta_{ij}$ and
\[
\partial_{\xi_i}(\xi_{a_1}\cdots\xi_{a_l})
= \sum_{k}(-1)^{k-1}\xi_{a_1}\cdots\xi_{a_{k-1}}
\partial_{\xi_i}(\xi_{a_k}) \xi_{a_{k+1}}\cdots\xi_{a_l}.
\]

The formulas below for type $A_{n-1}$ case can be obtained from
Lemmas~\ref{Bcomm:[eta,xi]}, \ref{Bcomm:[eta,f]}, and
Theorem~\ref{BDunkl:eta} with the removal of those terms involving
$\overline{s}_{ij}, \overline{s}_{ki},$ and the parameter $v$
therein.

\begin{Lemma}
Let $W=W_{A_{n-1}}$. Then the following holds in $\HO_W$ for $l \in \Z_+$ and $i\neq j$:
\begin{gather*}
\lbrack \eta_{i},\xi_{j}^l]_{+}
 =  \frac{u}{\xi_i^2-\xi_j^2} \big(\xi_i^{l+1}-\xi_j\xi_i^{l}
-\xi_i\xi_j^{l}+(-1)^l\xi_j^{l+1} \big)  s_{ij}, \\
\lbrack \eta_{i},\xi_{i}^l]_{+}
 = t\frac{\xi_{i}^l - (-\xi_{i})^l}{2\xi_i} + u\sum_{k\neq i}
\frac{1}{\xi_i^2-\xi_k^2}\big(\xi_i\xi_k^l -
\xi_k^{l+1}-(-1)^l\xi_i^{l+1} + \xi_k \xi_i^l \big) s_{ik}.
\end{gather*}
\end{Lemma}

\begin{Lemma}
Let $W=W_{A_{n-1}}$, and $f \in \Cl[\xi_1,\ldots,\xi_n]$. Then the
following identity holds in $\HO_W$:
\begin{gather*}
\lbrack \eta_{i},f]_{+}  =  t\frac{f-f^{\tau_i}}{2\xi_i}
+u \sum_{k\neq i} \frac{1}{\xi_i^2-\xi_k^2}\big(\big(\xi_i
- \xi_k\big)f^{s_{ik}} - \big(\xi_i f^{\tau_i} - \xi_k f^{\tau_k}\big)\big)s_{k i} .
\end{gather*}
\end{Lemma}

\begin{Theorem} \label{ADunkl:eta}
Let $W=W_{A_{n-1}}$. The action of $\eta_i$ on
$\Cl[\xi_1,\ldots,\xi_n]$ is realized as Dunkl operators as
follows:
\begin{gather*}
\eta_i  =  t\partial_{\xi_i}+ u\sum_{k\neq
i}\frac{1}{\xi_i^2-\xi_k^2}\big(\big(\xi_i -
\xi_k\big)s_{ik}-\big(\xi_i\tau_i - \xi_k \tau_k\big)\big).
\end{gather*}
\end{Theorem}

\subsubsection[Dunkl operator for type $B_n$ case]{Dunkl operator for type $\boldsymbol{B_n}$ case}

The proofs of Lemma~\ref{Bcomm:[eta,xi]} and \ref{Bcomm:[eta,f]}
are given in the Appendix.

\begin{Lemma} \label{Bcomm:[eta,xi]}
Let $W=W_{B_{n}}$. Then the following holds in $\HO_W$ for $l \in \Z_+$ and $i\neq j$:
\begin{gather*}
\lbrack \eta_{i},\xi_{j}^l]_{+}
 =  u\big(\xi_i^{l-1} -\xi_j \xi_i^{l-2} + \dots
+ (-1)^{l-1}\xi_j^{l-1}\big)\big(s_{ij}+\overline{s}_{ij}\big)\\
\phantom{\lbrack \eta_{i},\xi_{j}^l]_{+}}{}  =  u\big(\frac{1}{\xi_i^2-\xi_j^2} \big(\xi_i^{l+1}-\xi_j\xi_i^{l}
-\xi_i\xi_j^{l}+(-1)^l\xi_j^{l+1} \big) \big)\big(s_{ij}+\overline{s}_{ij}\big),\\
%
\lbrack \eta_{i},\xi_{i}^l]_{+}  =  t\frac{\xi_{i}^l -
(-\xi_{i})^l}{2\xi_i} + v \frac{\xi_{i}^l - (-\xi_{i})^l}{2\xi_i}
\tau_i\\
\phantom{\lbrack \eta_{i},\xi_{i}^l]_{+}=}{}
+ u \sum_{k\neq i}\big(\xi_k^{l-1} -\xi_i \xi_k^{l-2}
+ \dots + (-1)^{l-1}\xi_i^{l-1} \big)\big(s_{ik}+\overline{s}_{ik}\big)\\
\phantom{\lbrack \eta_{i},\xi_{i}^l]_{+}}{}  =  t\frac{\xi_{i}^l - (-\xi_{i})^l)}{2\xi_i}
+ v \frac{\xi_{i}^l - (-\xi_{i})^l}{2\xi_i} \tau_i \\
\phantom{\lbrack \eta_{i},\xi_{i}^l]_{+}=}{} +u\sum_{k\neq i}\frac{1}{\xi_i^2-\xi_k^2}\big(\xi_i\xi_k^l
- \xi_k^{l+1}-(-1)^l\xi_i^{l+1} + \xi_k \xi_i^l \big)\big(s_{ik}+\overline{s}_{ik}\big).
\end{gather*}
\end{Lemma}

\begin{Lemma} \label{Bcomm:[eta,f]}
Let $W=W_{B_{n}}$, and $f \in \Cl[\xi_1,\ldots,\xi_n]$. Then the
following identity holds in $\HO_W$:
\begin{gather*}
\lbrack \eta_{i},f]_{+}  =  t \frac{f-f^{\tau_i}}{2\xi_i}+ v
\frac{f-f^{\tau_i}}{2\xi_i} \tau_i\\
\phantom{\lbrack \eta_{i},f]_{+}  =}{}
+u \sum_{k\neq i}
\frac{1}{\xi_i^2-\xi_k^2}\big(\big(\xi_i - \xi_k\big)f^{s_{ik}} -
\big(\xi_i f^{\tau_i} - \xi_k
f^{\tau_k}\big)\big)\big(s_{ik}+\overline{s}_{ik}\big) .
\end{gather*}
\end{Lemma}

\begin{Theorem} \label{BDunkl:eta}
Let $W=W_{B_{n}}$. The action of $\eta_i$ on
$\Cl[\xi_1,\ldots,\xi_n]$ is realized as operators as follows:
\begin{gather*}
\eta_i  =  t\partial_{\xi_i}+ v\frac{1-\tau_i}{2\xi_i} +
u\sum_{k\neq i}\frac{2}{\xi_i^2-\xi_k^2}\big(\big(\xi_i -
\xi_k\big)s_{ik}-\big(\xi_i\tau_i - \xi_k \tau_k\big)\big).
\end{gather*}
\end{Theorem}

\begin{proof}
It suf\/f\/ices to check the formula for every monomial $f$. Consider
$f= \xi_1^{a_1}\cdots \xi_n^{a_n}$ where $a_i \in \Z_+$, and observe
that
\[
\partial_{\xi_i}(f)  = \frac{f-f^{\tau_i}}{2\xi_i},\qquad
\eta_i \cdot f = [\eta_i, f]_+ + (-1)^{a_1+\cdots+a_n}  f\cdot
\eta_i = [\eta_i,f]_+.
\]
The theorem now follows by Lemma~\ref{Bcomm:[eta,f]}.
\end{proof}

\subsubsection[Dunkl operator for type $D_n$ case]{Dunkl operator for type $\boldsymbol{D_n}$ case}

The formula below for the Dunkl operator type $D_n$ case is
obtained from their type $B_n$ counterparts (see
Theorem~\ref{BDunkl:eta}) by dropping the terms involving the
parameter $v$.

\begin{Theorem} \label{DDunkl:eta}
Let $W=W_{D_{n}}$. The action of $\eta_i$ on
$\Cl[\xi_1,\ldots,\xi_n]$ is realized as Dunkl operators as
follows:
\[
\eta_i = t \partial_{\xi_i}+ u \sum_{k\neq
i}\frac{2}{\xi_i^2-\xi_k^2}\big(\big(\xi_i -
\xi_k\big)s_{ik}-\big(\xi_i\tau_i - \xi_k \tau_k\big)\big).
\]
\end{Theorem}

\begin{Remark}\rm
Let $W=W_{A_{n-1}}, W_{B_n},$ or $W_{D_n}$. The Dunkl operators
$\eta_i$ anti-commute, i.e.\ $\eta_i \eta_j = -\eta_j \eta_i$
$(i\neq j)$. It is not easy to check this directly.
\end{Remark}
\subsection[An affine Hecke subalgebra]{An af\/f\/ine Hecke subalgebra}

In this subsection, we will show that the odd DaHa of type $A$
contains as a subalgebra   the degenerate af\/f\/ine Hecke algebra of
type $A$ introduced by Drinfeld and Lusztig \cite{Dr, Lu}. Let
\[
\mathfrak z_i = - \xi_i \eta_i + u\sum_{k<i}
s_{ki}.
\]

\begin{Lemma}
We have $[\mathfrak z_i, \mathfrak z_j] =0$, $\forall\, i, j.$
\end{Lemma}

\begin{proof}
Let us assume $i<j$. Then,
\begin{gather*}
 [\mathfrak z_i, \mathfrak z_j]
  =  \Bigg[- \xi_i\eta_i, - \xi_j \eta_j + u\sum_{k<j}
s_{kj}\Bigg]
 =  (\xi_i [\eta_i, \xi_j]_+ \eta_j -\xi_j [\eta_j, \xi_i]_+ \eta_i)
 - u[\xi_i\eta_i, s_{ij}]
 \\
 \phantom{[\mathfrak z_i, \mathfrak z_j]}{} =  u(\xi_i s_{ij} \eta_j -\xi_j s_{ij} \eta_i)
 - u (\xi_i\eta_i - \xi_j\eta_j) s_{ij}= 0.\tag*{\qed}
\end{gather*}\renewcommand{\qed}{}
\end{proof}

\begin{Lemma}\label{Hecke:identities}
The following identities hold:
    \begin{gather*}
        s_i \mathfrak z_i   = \mathfrak z_{i+1} s_i - u,\qquad
        s_i \mathfrak z_j   = \mathfrak z_j s_i \qquad (j \neq i, i+1).
    \end{gather*}
\end{Lemma}
\begin{proof}
Recall that $L_i :=\sum\limits_{k<i}s_{ki}$ is the Jucys--Murphy element,
and it is known that $s_i L_i = L_{i+1}s_i -1$ and $s_i L_j = L_j
s_i$ for $j\neq i, i+1$. The lemma follows from these relations.
\end{proof}

\begin{Proposition}
The $\mathfrak z_i$ $(1\le i \le n)$ and $S_n$ generate the degenerate affine Hecke
algebra.
\end{Proposition}
\begin{proof}
    The proposition follows from Theorem~\ref{PBW:ADBoddDaHa} and Lemma~\ref{Hecke:identities}.
\end{proof}

\appendix
\section{Appendix: proofs of several lemmas}
\label{sec:Appendix}
\subsection{Proofs of Lemmas in Section~\ref{sec:DaHa:doubleClifford}}

\subsubsection{Proof of Lemma~\ref{conj-inv-c}}

We will show that the relations (\ref{Byixj}) and (\ref{Byixi})
are invariant under the conjugation by ele\-ments~$c_l$ and~$e_l$,
$1\le l \le n$. We will only verify for the $c_l$ and leave the
similar verif\/ication for the $e_l$ to the reader. Also, the
verif\/ications for the invariants in type $A$ and $D$ under the
conjugation by $c_l$ and $e_l$ are similar and will be omitted.

Consider the relation (\ref{Byixj}) f\/irst. Clearly, (\ref{Byixj})
is invariant under the conjugation by $c_l$, and $e_l$ if $l\neq
i,j$. Moreover, we calculate that
\begin{gather*}
c_i   \text{(r.h.s. of (\ref{Byixj}))} c_i
  = u \big((1
+c_ic_j)(1 +e_je_i)s_{ij}
- (1-c_ic_j)(1-e_je_i)\overline{s}_{ij}\big) \\
\phantom{c_i   \text{(r.h.s. of (\ref{Byixj}))} c_i}{}
=  \lbrack y_{i}, x_{j}]
 = c_i \text{(l.h.s. of (\ref{Byixj}))} c_i, \\
c_j   \text{(r.h.s. of (\ref{Byixj}))} c_j
=  u\big((c_j c_i-1)(1+e_je_i)s_{ji}-(-c_{j} c_{i}-1)(1-e_je_i)\overline{s}_{ij}\big)\\
\phantom{c_j   \text{(r.h.s. of (\ref{Byixj}))} c_j}{}
=  -\lbrack y_{i},x_{j}]
 = c_j \text{(l.h.s. of (\ref{Byixj}))} c_j.
\end{gather*}
Thus, (\ref{Byixj}) is conjugation-invariant by all $c_l$.

Next, we will show that the relation (\ref{Byixi}) is invariant
under the conjugation by each $c_l$. Indeed, we have
\begin{gather*}
c_i   \text{(r.h.s. of  (\ref{Byixi}))} c_i \\
\qquad{} =  t e_i c_i - v c_i\tau_{i} c_i
-u\sum_{k\neq
i}c_i((1+c_{k}c_{i})(1+e_k e_i)s_{ki} +(1-c_{k}
c_{i})(1-e_k e_i)\overline{s}_{ki})c_i \\
\qquad{} = -tc_i e_i + v\tau_{i}
 -u\sum_{k\neq i}((c_i c_k -
1)(1+e_k e_i)s_{ki}+(-c_i c_k -1)(1-e_k e_i)\overline{s}_{ki})\\
\qquad{}= -tc_i e_i + v \tau_{i}
+ u\sum_{k\neq i}((1+ c_{k}c_{i})(1+e_k e_i)s_{ki}+(1-c_{k}%
c_{i})(1-e_k e_i)\overline{s}_{ki})\\
\qquad{}= -\lbrack y_{i},x_{i}]
 = c_i \text{(l.h.s. of (\ref{Byixi}))}
c_i.
\end{gather*}
For $j \neq i$, we have
\begin{gather*}
 c_j \text{(r.h.s. of (\ref{Byixi}))} c_j \\
\qquad{}  = tc_i e_i -  v \tau_{i} -u c_j((1+c_jc_i)(1+e_j e_i)s_{ji}
+(1-c_{j}c_{i})(1-e_j e_i)\overline{s}_{ji})c_j \\
\qquad\quad {}-u\sum_{k\neq i,j}c_j((1+c_{k}c_{i})(1+e_k e_i)s_{ki}+(1-c_{k}%
c_{i})(1-e_k e_i)\overline{s}_{ki})c_j\\
\qquad{}= tc_i e_i - v \tau_{i}
 -u ((c_jc_i+1)(1+e_j e_i)s_{ji}
+(-c_{j}c_{i}+1)(1-e_j e_i)\overline{s}_{ji})c_j \\
\qquad\quad {}-u\sum_{k\neq i,j}((1+c_{k}c_{i})(1+e_k e_i)s_{ki}+(1-c_{k}%
c_{i})(1-e_k e_i)\overline{s}_{ki})\\
\qquad{}= c_j \text{(l.h.s. of (\ref{Byixi}))} c_j.
\end{gather*}
Therefore, the lemma is proved.

\subsubsection{{Proof of Lemma~\ref{conj-inv-W}}}

We will show below that the relations (\ref{Byixj})--(\ref{Byixi})
are invariant under the conjugation by elements in $W_{B_n}$. The
proof can be readily modif\/ied to yield the Weyl group invariance
of the relations (\ref{Ayixj})--(\ref{Ayixi}) and
(\ref{Dyixj})--(\ref{Dyixi}) in type $A$ and $D$ cases respectively,
and we leave the details to the reader.

(i) We check the invariance of (\ref{Byixj}) under $W_{B_n}$.

Consider f\/irst the conjugation invariance by the transposition
$s_{l k}$. If $\{l, k\}\cap\{i,j\} = \varnothing$, then we have
\begin{gather*}
s_{l k}   \text{(r.h.s. of (\ref{Byixj}))} s_{l k}
 =  u\big((1+c_i
c_j)(1+e_j e_i)s_{ij}-(1-c_i c_j)(1-e_j e_i)\overline{s}_{ij}\big) \\
\phantom{s_{l k}   \text{(r.h.s. of (\ref{Byixj}))} s_{l k}}{}  =  \lbrack y_{i},x_{j}]
 = s_{l k} \text{(l.h.s. of (\ref{Byixj}))} s_{l k}.
\end{gather*}

If $\{l, k\}\cap\{i,j\} =\{j\}$, then we may assume $l = j$ and we
have
\begin{gather*}
s_{j k}  \text{(r.h.s. of (\ref{Byixj}))} s_{j k}
=  u\big(
(1+c_{i}c_{k})(1+e_k e_i)s_{ik}-(1-c_{i}c_{k})(1-e_k e_i)\overline{s}_{ik}\big)\\
\phantom{s_{j k}  \text{(r.h.s. of (\ref{Byixj}))} s_{j k}}{} =  \lbrack y_{i},x_{k}]
 = s_{j k} \text{(l.h.s. of (\ref{Byixj}))} s_{j k}.
\end{gather*}
We leave an entirely analogous computation when $\{l,
k\}\cap\{i,j\} =\{i\}$ to the reader.

Now, if $\{l, k\}=\{i,j\}$, then
\begin{gather*}
s_{i j}  \text{(r.h.s. of (\ref{Byixj}))} s_{i j}
 = u\big(
(1+c_j c_i)(1+e_{i}e_{j})s_{ij}-(1-c_j c_i)(1-e_{i}e_{j})\overline{s}_{ij}\big) \\
\phantom{s_{i j}  \text{(r.h.s. of (\ref{Byixj}))} s_{i j}}{} =  \lbrack y_{j},x_{i}]
 = s_{i j} \text{(l.h.s. of (\ref{Byixj}))} s_{i j}.
\end{gather*}
So (\ref{Byixj}) is invariant under the conjugation by each
transposition $s_{l k}$.

It remains to show that (\ref{Byixj}) is invariant under the
conjugation by the simple ref\/lection $s_n= \tau_n$. Observe that
(\ref{Byixj}) is clearly invariant under the conjugation by $s_n$
for $n \neq i,j$. Moreover, if $j=n$ then we have
\begin{gather*}
s_{n}   \text{(r.h.s. of (\ref{Byixj}))} s_{n}
   = u\big((1-c_i
c_j)(1-e_je_i)\overline{s}_{ij}-(1+c_i c_j)(1+e_je_i)s_{ij}\big)\\
\phantom{s_{n}   \text{(r.h.s. of (\ref{Byixj}))} s_{n}}{}
= -\lbrack y_{i},x_{j}]
 = s_{n} \text{(l.h.s. of (\ref{Byixj}))} s_{n}.
\end{gather*}
If $i=n$, then we have
\begin{gather*}
s_{n}   \text{(r.h.s. of (\ref{Byixj}))} s_{n}
 = u\big((1-c_i
c_j)(1-e_j e_i)\overline{s}_{ji}-(1+c_i c_j)(1+e_j e_i)s_{ij}\big)\\
\phantom{s_{n}   \text{(r.h.s. of (\ref{Byixj}))} s_{n}}{}
=  -\lbrack y_{i},x_{j}] = s_{n} \text{(l.h.s. of (\ref{Byixj}))}
s_{n}.
\end{gather*}
This completes (i).

 \vspace{.4cm}

(ii) We check the invariance of (\ref{Byixi}) under $W_{B_n}$.

Consider f\/irst the conjugation invariance  by  $s_{j l}$. If $\{j,
l\}\cap\{i\} = \varnothing$, then we have
\begin{gather*}
s_{j l}   (\text{r.h.s. of } (\ref{Byixi})) s_{j l} \\
\qquad{}  = t c_i e_i -  v \tau_{i}
-u\sum_{k\neq i,j,l}s_{j l}((1+c_{k}c_{i})(1+e_{k}e_{i})s_{ki}+(1-c_{k}%
c_{i})(1-e_{k}e_{i})\overline{s}_{ki})s_{j l}\\
\qquad\quad{}-u s_{j l}\left((1+c_j c_i)(1+e_j e_i)s_{j i}+(1-c_j c_i)(1-e_j e_i)
\overline{s}_{ji}\right)s_{j l} \\
\qquad\quad{} -us_{j l}\left((1+c_l c_i)(1+e_l e_i)s_{l i}+(1-c_l c_i)(1-e_l e_i)
\overline{s}_{li}\right)s_{j l}\\
\qquad{}= \lbrack y_{i},x_{i}]
=s_{j l} \text{(l.h.s. of (\ref{Byixi}))}s_{j l}.
\end{gather*}
If $\{j, l\}\cap\{i\} = \{i\}$, we may assume that $j = i$, and then
we have
\begin{gather*}
s_{i l}   \text{(r.h.s. of } (\ref{Byixi})) s_{i l}\\
\qquad{}=t c_l e_l-u s_{i l}\left((1+c_l c_i)(1+e_l e_i)s_{j i}
+(1-c_l c_i)(1-e_l e_i)\overline{s}_{l i}\right)s_{i l} \\
\qquad\quad{} -u\sum_{k\neq i,l} s_{i l}((1+c_{k}c_{l})(1+e_{k}e_{i})s_{kl}+(1-c_{k}%
c_{l})(1-e_{k}e_{i})\overline{s}_{kl})s_{i l} - v \tau_{l}\\
\qquad{}= \lbrack y_{l},x_{l}] =s_{i l} \text{(l.h.s. of (\ref{Byixi}))}s_{i l}.
\end{gather*}

It remains to show that (\ref{Byixi}) is invariant under the
conjugation by the simple ref\/lection $s_n \equiv \tau_n \in
W_{B_n}$. If $i \neq n$, we have
\begin{gather*}
s_{n}  \text{(r.h.s. of } (\ref{Byixi})) s_{n} \\
\qquad{}
= t c_i e_i-  v \tau_{i}
-u s_n\left((1+c_n c_i)(1+e_n e_i)s_{n
i}+(1-c_n
c_i)(1-e_n e_i)\overline{s}_{n i})\right)s_{n} \\
\qquad\quad{} -u\sum_{k\neq i,n}s_n ((1+c_{k}c_{i})(1+e_k e_i)s_{ki}+(1-c_{k}%
c_{i})(1-e_k e_i)\overline{s}_{ki})s_n \\
\qquad{}= - v \tau_{i}
 -u \left((1-c_n c_i)\overline{s}_{n
i}+(1+c_n c_i)s_{n i}
)\right)
 -u\sum_{k\neq i,n}((1+c_{k}c_{i})s_{ki}+(1-c_{k}%
c_{i})\overline{s}_{ki})\\
\qquad{}= \lbrack y_{i},x_{i}]
= s_{n} \text{(l.h.s. of (\ref{Byixi}))} s_{n}.
\end{gather*}
 If $i = n$, then
\begin{gather*}
s_{n}  \text{(r.h.s. of }  (\ref{Byixi})) s_{n} \\
\qquad{}= t c_n e_n-  v \tau_{n}
-u\sum_{k\neq n}((1-c_k c_n)(1-e_k e_n)\overline{s}_{kn }+(1+c_k c_n)(1+e_k e_n)s_{kn}) \\
\qquad{}= \lbrack y_{n},x_{n}]
= s_{n} \text{(l.h.s. of (\ref{Byixi}))} s_{n}.
\end{gather*}
This completes the proof of (ii). Hence the lemma is proved.

\subsubsection{Proof of Lemma~\ref{Jacobi}}

We will establish the Jacobi identity for $W = W_{B_n}$. The proof
can be easily modif\/ied for the cases of type $A$ and $D$, and we
leave the details to the reader.

The Jacobi identity trivially holds among triple $x_i$'s or triple
$y_i$'s.

Now, we consider the triple with two $y$'s and one $x$. The case
with two identical $y_i$ is trivial. So we f\/irst consider $x_i$,
$y_j$, and $y_l$ where $i,j,l$ are all distinct. The Jacobi
identity holds in this case since
\begin{gather*}
[x_i,[y_j,y_l]]  +[y_l,[x_i,y_j]]+[y_j,[y_l,x_i]] \\
\qquad{}= 0+ [y_l, -u\big((1+c_{j}c_{i})(1+e_i e_j)s_{ji}-(1-c_{j}
c_{i})(1-e_i e_j)\overline{s}_{ij}\big)] \\
\quad\qquad {}+ [y_j, u\left((1+c_{l}c_{i})(1-e_i e_l)s_{li}-(1-c_{l}
c_{i})(1-e_i e_l)\overline{s}_{il}\right)] =0.
\end{gather*}

Now for $i\neq j$, we have
\begin{gather*}
[x_i,[y_i,y_j ]] +[y_j,[x_i,y_i]]+[y_i,[y_j,x_i]]\\
\qquad{} =
  \Bigg[y_j, -tc_i e_i + u\sum_{k\neq i} \big ( (1+c_{k}c_{i})(1+e_k
e_i)s_{ki}+(1-c_{k}c_{i})(1-e_k e_i)
\overline{s}_{ki} \big )+  v\tau_{i}\Bigg]\\
 \qquad \quad{}+ [y_i,u\left((1+c_{j}c_{i})(1+e_i e_j)s_{ij}-(1-c_{j}
    c_{i})(1-e_i e_j)\overline{s}_{ij}\right)]\\
\qquad{}= \Bigg[y_j,u\sum_{k\neq i,j}\big ( (1+c_{k}c_{i})(1+e_{k}e_{i})s_{ki}
  +(1-c_{k}c_{i})(1-e_{k}e_{i})\overline{s}_{ki} \big ) \Bigg]\\
\qquad \quad{}+ [y_j,u \big ( (1+c_{j}c_{i})(1+e_j e_i)s_{ji}
  +(1-c_jc_{i})(1-e_je_i)\overline{s}_{ji} \big )]\\
 \qquad \quad{}+ [y_i,u\left( (1+c_{j}c_{i})(1+e_i e_j)s_{ij}-(1-c_{j}
    c_{i})(1-e_i e_j)\overline{s}_{ij}\right)]\\
\qquad{}= 0+ u \big (y_j(1+c_j c_{i})(1+e_je_i) s_{ji}+y_j(1-c_j
    c_{i})(1-e_je_i)\overline{s}_{ji} \big) \\%
\quad\qquad {}-u \big ((1+c_jc_{i})(1+e_je_i)s_{ji}y_j+(1-c_j
    c_{i})(1-e_je_i)\overline{s}_{ji}y_j \big ) \\
\quad \qquad {}+u\left( y_i(1+c_{j}c_{i})(1+e_ie_j) s_{ji}-y_i(1-c_{j}
    c_{i})(1-e_ie_j) \overline{s}_{ij}\right) \\%
\qquad \quad{}-u\big((1+c_{j}c_{i})(1+e_ie_j)s_{ji}y_i-(1-c_{j}
    c_{i})(1-e_ie_j)\overline{s}_{ij}y_i\big)
    = 0.
\end{gather*}
Thanks to the automorphism $\varpi$ of $\bH_W$ which switches
$x_i$ and $y_i$, we obtain the Jacobi identity with one $y$ and
two $x$'s from the above calculation. This completes the proof of
Lemma~\ref{Jacobi}.

\subsubsection{{Proof of Lemma~\ref{Bcomm:[y,x^l]}}}

We will proceed by induction on $l$. For $l=1$, then the equations
hold by (\ref{Byixj}) and (\ref{Byixi}). Now assume that the
statement is true for $l$. Then
    \begin{gather*}
    \lbrack y_{i},x_{j}^{l+1}]  =  \lbrack y_{i},x_{j}^{l}] x_j + x_j^l \lbrack y_{i},x_{j}]\\
\phantom{\lbrack y_{i},x_{j}^{l+1}]}{} =  u \left(\frac{x_j^l - x_i^l}{x_j - x_i} +
    \frac{x_j^l - (-x_i)^l}{x_j + x_i}c_i c_j\right)(1-e_i e_j)s_{i j}x_j\\
\phantom{\lbrack y_{i},x_{j}^{l+1}]=}{} - u \left(\frac{x_j^l - (-x_i)^l}{x_j + x_i} - \frac{x_j^l -
    x_i^l}{x_j - x_i}c_i c_j\right)(1+e_i e_j)\overline{s}_{i j}x_j\\
\phantom{\lbrack y_{i},x_{j}^{l+1}]=}{} + x_j^l u\big((1+c_{i}c_{j})(1+e_j e_i)s_{ij}-(1-c_{i}c_{j})(1-e_j e_i)\overline{s}_{ij}\big)\\
\phantom{\lbrack y_{i},x_{j}^{l+1}]}{}= u \left(\frac{x_j^{l+1} - x_i^{l+1}}{x_j - x_i} +
    \frac{x_j^{l+1} - (-x_i)^{l+1}}{x_j + x_i}c_i c_j\right)(1-e_i e_j)s_{i j}\\
\phantom{\lbrack y_{i},x_{j}^{l+1}]=}{} - u \left(\frac{x_j^{l+1} - (-x_i)^{l+1}}{x_j + x_i}-
    \frac{x_j^{l+1} - x_i^{l+1}}{x_j - x_i}c_i c_j\right)(1+e_i e_j)\overline{s}_{i j},
\\
    \lbrack y_{i},x_{i}^{l+1}]  =  \lbrack y_{i},x_{i}^{l}] x_i + x_i^l \lbrack y_{i},x_{i}]\\
\phantom{\lbrack y_{i},x_{i}^{l+1}]}{} =
 t c_i e_i \frac{x_i^l - (x_i^l)^{\tau_i}}{2} - v\frac{x_i^l - (x_i^l)^{\tau_i}}{2} \tau_i\\
\phantom{\lbrack y_{i},x_{i}^{l+1}]=}{}-u \sum_{k\neq i} \left(\frac{x_i^l - x_k^l}{x_i - x_k} +
    \frac{x_i^l - (-x_k)^l}{x_i + x_k}c_k c_i\right)(1+e_k e_i)s_{k i}x_i \\
\phantom{\lbrack y_{i},x_{i}^{l+1}]=}{} - u \sum_{k\neq i} \left(\frac{x_i^l - (-x_k)^l}{x_i + x_k} -
    \frac{x_i^l - x_k^l}{x_i - x_k}c_k c_i\right)(1-e_k e_i)\overline{s}_{k i}x_i \\
\phantom{\lbrack y_{i},x_{i}^{l+1}]=}{} -u x_i^l\sum_{k\neq i}((1+c_{k}c_{i})(1+e_k e_i)s_{ki}+(1-c_{k}
    c_{i})(1-e_k e_i)\overline{s}_{ki})
     + t x_i^l c_i e_i - v x_i^l \tau_i \\
\phantom{\lbrack y_{i},x_{i}^{l+1}]}{} =  t c_i e_i \frac{x_i^{l+1} - (x_i^{l+1})^{\tau_i}}{2x_i}
    - v\frac{x_i^{l+1} - (x_i^{l+1})^{\tau_i}}{2x_i} \tau_i\\
\phantom{\lbrack y_{i},x_{i}^{l+1}]=}{} -u \sum_{k\neq i} \left(\frac{x_i^{l+1} - x_k^{l+1}}{x_i - x_k} +
    \frac{x_i^{l+1} - (-x_k)^{l+1}}{x_i + x_k}c_k c_i\right)(1+e_k e_i)s_{k i} \\
\phantom{\lbrack y_{i},x_{i}^{l+1}]=}{} - u \sum_{k\neq i} \left(\frac{x_i^{l+1} - (-x_k)^{l+1}}{x_i + x_k}
    - \frac{x_i^{l+1} - x_k^{l+1}}{x_i - x_k}c_k c_i\right)(1-e_k e_i)\overline{s}_{k i}.
\end{gather*}
This completes the proof.

\subsubsection{Proof of Lemma~\ref{Bcomm:[y,f]}}

It suf\/f\/ices to check the formula for every monomial $f$. First, we
consider the monomial $g = \prod\limits_{j\neq i} x_j^{a_j}$. By
induction and Lemma~\ref{Bcomm:[y,x^l]}, we can show that the
formula holds for the monomial of the form $g = \prod\limits_{j\neq i}
x_j^{a_j}$ (the detail of the induction step does not dif\/fer much
from the following calculation). Now consider the monomial $f
=x_i^{l} g$.
\begin{gather*}
\lbrack y_i,f]  =  \lbrack y_i,x_i^{l}]g + x_i^l \lbrack y_i, g]
\\
\phantom{\lbrack y_i,f]}{} =  t c_i e_i \frac{x_i^l-(-x_i)^l}{2x_i}g - v \frac{x_i^l - (-x_i)^l}{2x_i}\tau_i g \\
\phantom{\lbrack y_i,f]=}{} -u \sum_{k\neq i} \left(\frac{x_i^l - x_k^l}{x_i - x_k} +
        \frac{x_i^l - (-x_k)^l}{x_i + x_k}c_k c_i\right)(1+e_k e_i)s_{k i}g \\
\phantom{\lbrack y_i,f]=}{} - u \sum_{k\neq i} \left(\frac{x_i^l - (-x_k)^l}{x_i + x_k} -
        \frac{x_i^l - x_k^l}{x_i - x_k}c_k c_i\right)(1-e_k e_i)\overline{s}_{k i}g \\
\phantom{\lbrack y_i,f]=}{} +tx_i^lc_i e_i \frac{g-g^{\tau_i}}{2x_i}- v x_i^l\frac{g-g^{\tau_i}}{2x_i}\tau_i \\
\phantom{\lbrack y_i,f]=}{} -u\sum_{k\neq i} x_i^l\left(\frac{g- g^{s_{ki}}}{x_i - x_k} +
    \frac{g - g^{\overline{s}_{k i}}}{x_i + x_k}c_k c_i\right)(1+e_k e_i)s_{k i} \\
\phantom{\lbrack y_i,f]=}{} - u \sum_{k\neq i} x_i^l\left(\frac{g - g^{\overline{s}_{k
    i}}}{x_i + x_k} - \frac{g - g^{s_{k i}}}{x_i - x_k}c_k
    c_i\right)(1-e_k e_i)\overline{s}_{k i}\\
\phantom{\lbrack y_i,f]}{} =  tc_i e_i \frac{f-f^{\tau_i}}{2x_i}- v \frac{f-f^{\tau_i}}{2x_i}\tau_i
    -u \sum_{k\neq i} \left(\frac{f - f^{s_{ki}}}{x_i - x_k} +
    \frac{f - f^{\overline{s}_{k i}}}{x_i + x_k}c_k c_i\right)(1+e_k e_i)s_{k i} \\
\phantom{\lbrack y_i,f]=}{} - u \sum_{k\neq i} \left(\frac{f - f^{\overline{s}_{k i}}}{x_i +
    x_k} - \frac{f - f^{s_{k i}}}{x_i - x_k}c_k
    c_i\right)(1-e_k e_i)\overline{s}_{k i}.
\end{gather*}
So the lemma is proved.

\subsection{Proofs of Lemmas in Section~\ref{sec:oddDaHa}}

\subsubsection{{Proof of Lemma~\ref{Bcomm:[eta,xi]}}}

We will proceed by induction on $l$. For $l=1$, then the equations
hold by the def\/inition of $\HO_W$. Now assume that the statement
is true for $l$. Then
    \begin{gather*}
    \lbrack \eta_{i},\xi_{j}^{l+1}]_{+}
     =  \lbrack \eta_{i},\xi_{j}^{l}]_{+} \xi_j + (-1)^{l}\xi_j^l \lbrack \eta_{i},\xi_{j}]_{+}\\
\phantom{\lbrack \eta_{i},\xi_{j}^{l+1}]_{+}}{} = u\frac{1}{\xi_i^2-\xi_j^2} \big(\xi_i^{l+1}-\xi_j\xi_i^{l}
        -\xi_i\xi_j^{l}+(-1)^l\xi_j^{l+1} \big) \big(s_{ij}+\overline{s}_{ij}\big)\xi_j\\
\phantom{\lbrack \eta_{i},\xi_{j}^{l+1}]_{+}=}{}+u \frac{ (-\xi_j)^l}{\xi_i^2-\xi_j^2}(\xi_i^2-\xi_j^2) \big(s_{ij}+\overline{s}_{ij}\big)\\
\phantom{\lbrack \eta_{i},\xi_{j}^{l+1}]_{+}}{} =  u\frac{1}{\xi_i^2-\xi_j^2} \big(\xi_i^{l+2}-\xi_j\xi_i^{l+1}
        -\xi_i\xi_j^{l+1}+(-1)^{l+1}\xi_j^{l+2} \big) \big(s_{ij}+\overline{s}_{ij}\big),
\\
        \lbrack \eta_{i},\xi_{i}^{l+1}]_{+}
 =  \lbrack \eta_{i},\xi_{i}^{l}]_{+} \xi_i + (-1)^{l}\xi_i^l \lbrack \eta_{i},\xi_{i}]_{+}\\
\phantom{\lbrack \eta_{i},\xi_{i}^{l+1}]_{+}}{} =  t\frac{\xi_{i}^l - (\xi_{i}^l)^{\tau_i}}{2\xi_i}\xi_i
+ v \frac{\xi_{i}^l - (\xi_{i}^l)^{\tau_i}}{2\xi_i} \tau_i \xi_i \\
\phantom{\lbrack \eta_{i},\xi_{i}^{l+1}]_{+}=}{}
+u\sum_{k\neq i}\frac{1}{\xi_i^2-\xi_k^2}
\big(\xi_i\xi_k^l - \xi_k^{l+1}-(-1)^l\xi_i^{l+1} + \xi_k \xi_i^l \big)
\big(s_{ik}+\overline{s}_{ik}\big)\xi_i\\
\phantom{\lbrack \eta_{i},\xi_{i}^{l+1}]_{+}=}{} + t (-\xi_i)^l + v (-\xi_i)^l \tau_i
+ u\sum_{k \ne i}\frac{(-\xi_i)^l }{\xi_i^2-\xi_k^2}\big(\xi_i^2-\xi_k^2\big)\big(s_{ki}
+ \overline{s}_{ij}\big)\\
\phantom{\lbrack \eta_{i},\xi_{i}^{l+1}]_{+}}{}
=  t\frac{\xi_{i}^{l+1} - (\xi_{i}^{l+1})^{\tau_i}}{2\xi_i}
+ v \frac{\xi_{i}^{l+1} - (\xi_{i}^{l+1})^{\tau_i}}{2\xi_i} \tau_i \\
\phantom{\lbrack \eta_{i},\xi_{i}^{l+1}]_{+}=}{}
 +u\sum_{k\neq i}\frac{1}{\xi_i^2-\xi_k^2}\big(\xi_i\xi_k^{l+1}
- \xi_k^{l+2}-(-1)^{l+1}\xi_i^{l+2} + \xi_k \xi_i^{l+1}
\big)\big(s_{ik}+\overline{s}_{ik}\big).
    \end{gather*}
This completes the proof.

\subsubsection{{Proof of Lemma~\ref{Bcomm:[eta,f]}}}
It suf\/f\/ices to check the formula for every monomial $f$. First, we
consider the monomial $g = \prod\limits_{j\neq i} \xi_j^{a_j}$. By
induction and Lemma \ref{Bcomm:[eta,xi]}, we can show that the
formula holds for the monomial of the form $g = \prod\limits_{j\neq i}
\xi_j^{a_j}$ (the detail of the induction step does not dif\/fer much
from the following calculation). Now consider the monomial $f
=\xi_i^{l} g$.
\begin{gather*}
\lbrack \eta_i,f]_{+}
  =  \lbrack
\eta_i,\xi_i^{l}]_{+}g +(-1)^l \xi_i^l \lbrack \eta_i, g]_{+}
\\
\phantom{\lbrack \eta_i,f]_{+}}{}
= t \frac{\xi_i^l g - (\xi_i^l g)^{\tau_i}}{2\xi_i}
+ v \frac{\xi_i^l g - (\xi_i^l g)^{\tau_i}}{2\xi_i} \tau_i \\
\phantom{\lbrack \eta_i,f]_{+}=}{}
+u\sum_{k\neq i}\frac{1}{\xi_i^2-\xi_k^2}\big(\xi_i\xi_k^l
- \xi_k^{l+1}-(-1)^l\xi_i^{l+1} + \xi_k \xi_i^l \big)\big(s_{ik}+\overline{s}_{ik}\big)g \\
\phantom{\lbrack \eta_i,f]_{+}=}{}
 +u \sum_{k\neq i} \frac{(-\xi_i)^l}{\xi_i^2-\xi_k^2}\big(\big(\xi_i
- \xi_k\big)g^{s_{ik}} - \big(\xi_i g^{\tau_i}
- \xi_k g^{\tau_k}\big)\big)\big(s_{ik}+\overline{s}_{ik}\big)\\
\phantom{\lbrack \eta_i,f]_{+}}{}
=  t \frac{\xi_i^l g - (\xi_i^l g)^{\tau_i}}{2\xi_i} + v \frac{\xi_i^l g
- (\xi_i^l g)^{\tau_i}}{2\xi_i} \tau_i\\
\phantom{\lbrack \eta_i,f]_{+}=}{}
+u\sum_{k\neq i}\frac{1}{\xi_i^2-\xi_k^2}
\big( (\xi_i -\xi_k)(\xi_i^l g)^{s_{ik}} -(\xi_i (\xi_i^l)^{\tau_i}
- \xi_k(\xi_i^l)^{\tau_k} )g^{s_{ik}} \big)\big(s_{ik}+\overline{s}_{ik}\big)\\
\phantom{\lbrack \eta_i,f]_{+}=}{}
 +u \sum_{k\neq i} \frac{(-\xi_i)^l}{\xi_i^2-\xi_k^2}\big(\big(\xi_i
- \xi_k\big)g^{s_{ik}} - \big(\xi_i g^{\tau_i} - \xi_k g^{\tau_k}\big)\big)
\big(s_{ik}+\overline{s}_{ik}\big)\\
\phantom{\lbrack \eta_i,f]_{+}}{}
=   t \frac{\xi_i^l g - (\xi_i^l g)^{\tau_i}}{2\xi_i}
+ v \frac{\xi_i^l g - (\xi_i^l g)^{\tau_i}}{2\xi_i} \tau_i\\
\phantom{\lbrack \eta_i,f]_{+}=}{}
 +u \sum_{k\neq i} \frac{1}{\xi_i^2-\xi_k^2}\big(\big(\xi_i
- \xi_k\big) (\xi_i^l g)^{s_{ik}} - \big(\xi_i (\xi_i^l g)^{\tau_i}
- \xi_k (\xi_i^l g)^{\tau_k}\big)\big)\big(s_{ik}+\overline{s}_{ik}\big).
\end{gather*}
So the lemma is proved.

\subsection*{Acknowledgements}
This research is partially supported by NSF grant DMS-0800280. The
main results of this paper for type $A$ were obtained at MSRI in
2006.

\pdfbookmark[1]{References}{ref}
\LastPageEnding

\end{document}